\RequirePackage{etoolbox}
\csdef{input@path}{{style/}{graphics/}}
\documentclass[aap,MSNbibl,seceqn,dvips]{arximspdf}
\makeatletter
   \@ifpackageloaded{graphicx}{}{\usepackage{graphicx}}
\makeatother

\doi{10.1214/14-AAP1028} 
\volume{25}
\issue{3}
\pubyear{2015}
\firstpage{1475}
\lastpage{1512}

\makeatletter

\newcommand{\eqref}[1]{(\ref{#1})}
\newcommand{\IR}{\mathbb{R}}

\newcommand{\B}{\mathcal{B}}

\newcommand{\F}{\mathcal{F}}

\newcommand{\ep}{\varepsilon}

\newcommand{\weak}{\mathrm{w}}

\newtheorem{stat}{Statement}[section]
\newproclaim{examp}[stat]{Example}
\newproclaim{assump}{Assumption}[section]

\newtheorem{prop}[stat]{Proposition}

\newtheorem{thmm}[stat]{Theorem}
\newtheorem{lemma}[stat]{Lemma}
\newproclaim{remark}[stat]{Remark}
\newproclaim{deff}[stat]{Definition} 
\makeatother

\begin{document}
\begin{frontmatter}

\title{A quickest detection problem with an~observation~cost}
\runtitle{A quickest detection problem}

\begin{aug}
\author[A]{\fnms{Robert C.} \snm{Dalang}\corref{}\thanksref{T2}\ead[label=e1]{robert.dalang@epfl.ch}}
\and
\author[B]{\fnms{Albert N.} \snm{Shiryaev}\thanksref{T3}\ead[label=e2]{albertsh@mi.ras.ru}}
\thankstext{T2}{Supported in part by the Swiss National Foundation for
Scientific Research.}
\thankstext{T3}{Supported in part by the Russian Foundation for Basic
Research (project 11-01-00949-a) and PreMoLab,
Moscow Institute of Physics and Technology (project 11.G34.31.0073).}
\runauthor{R. C. Dalang and A. N. Shiryaev}
\affiliation{Ecole Polytechnique F\'ed\'erale de Lausanne and\break Steklov
Mathematical Institute, Moscow}
\address[A]{Institut de math\'ematiques\\
Ecole Polytechnique F\'ed\'erale de Lausanne\\
Station 8\\
CH-1015 Lausanne\\
Switzerland\\
\printead{e1}}
\address[B]{Department of Probability Theory\\
\quad and Mathematical Statistics\\
Steklov Mathematical Institute\\
8, Gubkina str.\\
119991 Moscow\\
Russia\\
\printead{e2}}
\end{aug}

\received{\smonth{1} \syear{2013}}
\revised{\smonth{10} \syear{2013}}

%
\begin{abstract}
In the classical quickest detection problem, one must detect as quickly
as possible when a Brownian motion without drift ``changes'' into a
Brownian motion with positive drift. The change occurs at an unknown
``disorder'' time with exponential distribution. There is a penalty for
declaring too early that the change has occurred, and a cost for late
detection proportional to the time between occurrence of the change and
the time when the change is declared. Here, we consider the case where
there is also a cost for observing the process. This stochastic control
problem can be formulated using either the notion of strong solution or
of weak solution of the s.d.e. that defines the observation process. We
show that the value function is the same in both cases, even though no
optimal strategy exists in the strong formulation. We determine the
optimal strategy in the weak formulation and show, using a form of the
``principle of smooth fit'' and under natural hypotheses on the
parameters of the problem, that the optimal strategy takes the form of
a two-threshold policy: observe only when the posterior probability
that the change has already occurred, given the observations, is larger
than a threshold $A\geq0$, and declare that the disorder time has
occurred when this posterior probability exceeds a threshold $B \geq
A$. The constants $A$ and $B$ are determined explicitly from the
parameters of the problem.
\end{abstract}

%

\end{frontmatter}

\section{Introduction}\label{sec1}

The classical quickest detection problem \cite{shiryaev}, Chapter~4.4,
is as follows. One observes a stochastic process $X=(X_t)_{t \geq0}$
that solves the stochastic differential equation (s.d.e.)
%
%
\begin{equation}
\label{e1.1} dX_t = r 1_{\{\theta\leq t\}} \,dt + \sigma
\,dW_t.
\end{equation}
Here, $r > 0$, $\sigma> 0$, $W=(W_t)_{t \geq0}$ is a standard
Brownian motion,
and $\theta$ is a nonnegative random variable that is independent of
$(W_t)$, sometimes called a ``disorder time,'' or a ``change point.'' The
random variable $\theta$ is not observed directly, but only through its
effect on the sample paths of $X$. When $t < \theta$, the observer is
simply watching a Brownian motion, but when $t \geq\theta$, a drift
(or signal) with intensity $r$ appears. The observer seeks to detect as
quickly as possible the appearance of this signal, while keeping
sufficiently low the probability of a ``false alarm,'' that is,
declaring that the signal has appeared when, in fact, it has not.
Typically, the distribution of $\theta$ is assumed known, and, given
$\theta>0$, even equal to an exponential distribution with known
parameter $\lambda> 0$; see \cite{PS} for many variations on this
problem and for numerous references.
%

In this paper, we consider the situation where there is an \emph
{observation cost} $b\geq0$ per unit time, and the observer can choose
to observe or not. When he does not observe, the process $X$ is
constant $(dX_t = 0)$, and when he does observe, $X$ satisfies~(\ref
{e1.1}). The objective is to detect the appearance of the signal as
quickly as possible, while keeping low the probability of false alarm
\emph{and} the cost of observation. Therefore, the problem is no longer
an optimal stopping problem but an \emph{optimal stopping/control problem}.

There are several papers in the literature that consider this type of
problem, in which there is either a cost or a constraint on
observations. A discussion already appears in Bather \cite{bather},
with a precise continuous-time formulation given in Balmer \cite
{balmer1,balmer2}: he allows only a restricted class of policies and
uses a different cost function than the one we define in \eqref{cost}
below. Dayanik \cite{dayanik} considers a continuous-time problem with
observations allowed only at fixed times. Banerjee and Veeravalli \cite
{BV1,BV2} consider a discrete-time formulation, in which observations
are costly only if they occur before the alarm time. They show that a
two-threshold policy is asymptotically optimal. Finally, Bayraktar and
Kravitz \cite{BK} consider a continuous-time problem in which
observations are allowed only at a discrete set of times that is
determined adaptively.

In this paper, we consider that the control $h=(h_t)_{t \geq0}$ is a
$[0,1]$-valued process, where $h_t = 1$ means that observation occurs,
and $h_t = 0$ means absence of observation. Therefore, the observation
process is described by the stochastic differential
%
%
\begin{equation}
\label{e1.2} dX_t = r h_t 1_{\{\theta\leq t \}} \,dt +
\sigma \sqrt{h_t} \,dW_t,\qquad X_0 = 0.
\end{equation}
Note that when $h_t \in\{0,1\}$, the square-root has no effect.
However, it will be convenient during the resolution of the problem to
consider also $h_t \in[0,1]$, and since we are free to decide the
formulation when $0 < h_t < 1$, we have chosen to use~(\ref{e1.2}).


We assume that all objects are defined on a filtered probability space
$(\Omega, \F,  (\F_t)_{t\ge0}, P)$. Therefore, $W_t=W_t(\omega)$,
$\theta
=\theta(\omega)$, and $h_t=h_t(\omega)$.
The assumption that $h_t$ depends on $\omega$ ($h_t=h_t(\omega)$) does
not create difficulties with definition of the (``$h$-controlled'')
process $X$ via formula \eqref{e1.2}. However, we must define precisely
what information the observer can use to decide to switch from one
value of $h_t(\omega)$ to another.

It is reasonable to assume that the control function $h_t$ depends on
$\omega$ via the observation process: $h_t(\omega)=h_t(X(\omega))$. In
this case, the s.d.e. \eqref{e1.2} will take the form
%
%
\begin{equation}
\label{e1.3-ASh} dX_t =rh_t(X)1_{\{\theta\le t\}} \,dt +\sigma
\sqrt{h_t(X)} \,dW_t,
\end{equation}
and, inevitably, we have to explain how to formulate this s.d.e. and
give a precise definition of the control $h=(h_t(X))_{t\ge0}$.

These questions are considered in Section~\ref{sec2}, where we give two precise
but distinct formulations of the notion of a solution of equation
\eqref
{e1.3-ASh}, according to whether we interpret $X$ as a \textit{strong} or
\textit{weak} solution of \eqref{e1.3-ASh}.
%
%
Then we derive some preliminary properties of the sufficient statistic
$\pi^h_t$, which is the conditional probability, given the observations
$(X_s, s \in[0, t])$, that $\theta\leq t$. In Section~\ref{sec3}, we
study the law of $\pi^h_t$, writing it, and the likelihood ratio
$\varphi^h_t = \pi^h_t /(1-\pi^h_t)$, as solutions of diffusion
equations in the filtration $\F^X$ of the observed process. In this
section, we also establish, in the spirit of \cite{elkaroui} and \cite
{OS}, a ``verification lemma'' (Lemma~\ref{lem8}) that gives sufficient
conditions for the optimality of a strategy.

In Section~\ref{sec4}, we give the form of a candidate optimal
strategy and associated candidate value function and derive the
ordinary differential equations with two free boundaries that
characterize this function. These are completed by imposing boundary
conditions that imply continuity and an appropriate degree of
smoothness at the boundaries; see \eqref{e1.10}--\eqref{e1.13a}. These
equations are then solved completely, up to the resolution of a
transcendental equation; see \eqref{e1.24}. The form of the solution
depends on the value of the observation cost $b$, and it turns out that
there are three regimes: if $b$ is large enough, then it is best never
to observe, and to stop simply when the posterior probability $\pi^h_t$
exceeds a certain threshold $B\in\,]0,1[$. For smaller positive values
of $b$, there are two thresholds $0< A<B<1$ such that it is best not to
observe when $\pi^h_t \leq A$, to observe when $\pi^h_t \in\,]A,B[$ and
to declare an alarm when $\pi^h_t \geq B$. The candidate value function
is given in Propositions \ref{lem10} and \ref{lem11}, depending on the
size of $b$. The third regime is when $b=0$, which is the classical
case of \cite{shiryaev} and corresponds to $0=A < B <1$.

For small positive values of $b$, the candidate value function and
optimal strategies are such that it is not clear whether an optimal
strategy does indeed exist! In fact, in the strong formulation, no
optimal strategy exists in general, but such an optimal strategy does
exist in the weak formulation. It turns out, however, that the value
function is the same in both formulations. We discuss this question at
the end of Section~\ref{sec4}.

In Section~\ref{sec5}, we show that the candidate value function of
Section~\ref{sec4} is indeed the value function in both the weak and
strong formulations (Theorems \ref{thmvf} and \ref{thmvf1}). However,
because of the absence of an optimal strategy in the strong
formulation, it is not possible to conclude directly from a
``verification lemma'' (Lemma~\ref{lem8}) that the candidate value
function is indeed the value function in the strong formulation.
Therefore, we use a different approach in Theorem~\ref{thmvf1}: for
$\ep>0$, we consider strategies that approximate the candidate optimal
strategy but are defined via s.d.e.'s with sufficiently smooth
coefficients. We then compute explicitly the cost associated with these
strategies. This requires computing the expected time to hit a
threshold, which, in turn, requires solving another o.d.e. [given in
\eqref{e5.72}]. We do this in Section~\ref{sec5}, and in
Proposition~\ref{lem15}, we show by direct calculation that the
expected costs of
the approximately optimal strategies converge to the candidate value
function, proving that this is indeed the value function in the strong
formulation.

\section{Stating the problem}\label{sec2}

Consider a filtered probability space $(\Omega, \F,\break  (\F_t)_{t\ge0},
P)$ with a filtration $(\F_t)$ (satisfying the \emph{usual hypotheses}
\cite{protter}). Let $\theta$ be a random variable defined on $\Omega$
that is $\F_0$-measurable. We assume that there are $\pi_0 \in[0,1]$
and $\lambda>0$ such that
%
%
\begin{equation}
\label{e1.2a0} P\{\theta= 0\} = \pi_0\quad \mbox{and}\quad P\{\theta> x \mid
\theta>0\} = e^{-\lambda x}.
\end{equation}

We let $W=(W_t)_{t\ge0}$ be a standard Brownian motion adapted to $(\F
_t)_{t\ge0}$ such that for all $t \geq0$, the process $(W_{s+t} -
W_t, s \geq0)$ is independent of $\F_t$. In particular, $(W_t)_{t\ge
0}$ is independent of $\theta$.

\subsection*{Controls and stopping times}

%
\begin{deff}\label{def21}
A progressively measurable process $h=(h_t(\omega))_{t\ge0}$
defined on $(\Omega, \F, (\F_t)_{t\ge0}, P)$ with values
in $[0,1]$ will be called a \emph{stochastic control}.
\end{deff}

Let $C(\IR_+, \IR)$ denote the space of continuous functions from
$\IR
_+$ to $\IR$.

%
\begin{deff}\label{def-2-2}
A \emph{canonical control} $h=(h_t(x))_{t\ge0}$ is a map $(t, x)
\mapsto h_t(x)$ from $\IR_+ \times C(\IR_+, \IR)$ to $[0,1]$ 
that is progressively measurable for the canonical filtration on $C(\IR
_+, \IR)$.

A \emph{canonical stopping time} $\tau=\tau(x)$ is a random variable
$\tau\colon C(\IR_+, \IR) \to\IR_+$ that is a stopping time relative
to the canonical filtration on $C(\IR_+, \IR)$.
%
\end{deff}

%
\begin{deff}\label{def-2-3}
A stochastic control $h=(h_t(\omega))_{t\ge0}$ is called an \emph
{admissible control} if it has the form $h_t(\omega) = h_t(X(\omega))$
for a canonical control $h_t(x)$, and the s.d.e.
%
%
\begin{equation}
\label{e1.2a} dX_t = r h_t(X) 1_{\{\theta\leq t\}} \,dt +
\sigma\sqrt{h_t(X)} \,dW_t,\qquad X_0 = 0,
\end{equation}
admits a strong solution in the sense of the next definition
(Definition~\ref{def-2-4}).
\end{deff}

%
\begin{deff}\label{def-2-4}
Assume that a filtered probability space $(\Omega, \F,\break (\F_t)_{t\ge0},
P)$ is given a priori together with a random variable $\theta=\theta
(\omega)$ which is $\F_0$-measurable and satisfies \eqref{e1.2a0}, %
and with a Brownian motion $W(\omega)=(W_t(\omega))_{t\ge0}$ such that
$W_t$ is $\F_t$-measurable, for all $t\ge0$.

A \emph{strong solution} of the s.d.e. \eqref{e1.2a} is a continuous
stochastic process $X=(X_t(\omega))_{t\ge0}$ that satisfies \eqref
{e1.2a} and $X_t$ is $\F_t$-measurable, for all $t\ge0$.
\end{deff}

One may consider also the case where \eqref{e1.2a} has a weak solution.

%
\begin{deff}\label{def-2-5}
We assume that a canonical control $h=(h_t(x))_{t\ge0}$ and the law of
$\theta$ in \eqref{e1.2a0} are given a priori. A \emph{weak solution}
of the s.d.e. \eqref{e1.2a} is a system of the following objects:
\begin{itemize}
\item[--] a filtered probability space $(\Omega, \F, (\F_t)_{t\ge0}, P)$
(which is not given a priori);

\item[--] a Brownian motion $W=(W_t)_{t\ge0}$ such that $W_t$ is $\F
_t$-measurable, for all $t\ge0$;

\item[--] an $\F_0$-measurable random variable $\theta$ with the law
specified in \eqref{e1.2a0};

\item[--] an $(\F_t)_{t\ge0}$-adapted process $X=(X_t)_{t\ge0}$ which
satisfies the s.d.e. \eqref{e1.2a}, that is, for all $t\ge0$
%
%
\begin{equation}
\label{e-2-3-ASh} X_t =\int_0^tr
h_s(X) 1_{\{\theta\leq s\}} \,ds +\int_0^t
\sigma\sqrt{h_s(X)} \,dW_s.
\end{equation}
\end{itemize}
\end{deff}

%
\begin{deff}\label{def-2-6}
For the case of strong solutions, a \emph{strategy} is a pair $(h,\tau
)$, where $h=(h_t(X(\omega)))_{t\ge0}$, $\tau=\tau(X(\omega))$ for some
canonical control $(h_t(x))_{t\ge0}$ and canonical stopping time $\tau(x)$.

For the case of weak solutions, $(h,\tau,X)$ is called a \emph
{control system}.
\end{deff}

\subsection*{Cost}
%
\begin{deff}\label{def-2-7}
The \emph{cost} associated with a strategy $(h, \tau)$ or a control
system $(h, \tau,X)$ is
%
%
\begin{eqnarray}
\label{cost} C(h, \tau) &=& C(h, \tau,X)
\nonumber
\\[-8pt]
\\[-8pt]
\nonumber
& =& 1_{\{\tau(X) < \theta\}} + a \bigl(\tau(X)-
\theta \bigr)1_{\{\tau(X) \geq\theta\}} + b \int^{\tau(X)}_0
h_t(X) \,dt,
\end{eqnarray}
where $a>0$, so as to penalize late detection of the alarm time $\theta
$, and $b \geq0$.
\end{deff}

Since the case $b=0$ is covered in \textup{\cite{shiryaev},
Chapter~4.4}, we will focus on the case $b>0$.

\subsection*{Objective}

Our first objective is to find the value
\[
V = \inf_{(h,\tau)} E \bigl(C(h, \tau) \bigr),
\]
where the infimum is over all strategies, and to find an \emph{optimal}
strategy $(h^*,\tau^*)$ that achieves this infimum, or at least, to
find a strategy that is within $\ep$ of this infimum ($\ep>0$). A
second objective is to find the value
\[
V^{\weak} = \inf_{(h,\tau,X)} E \bigl(C(h, \tau,X) \bigr),
\]
where the infimum is over all control systems, and an optimal control
system $(h^*,\tau^*, X^*)$. Clearly, $V^{\weak} \leq V$.

\subsection*{Dependence on $\pi_0$}

The quantities $V$ and $V^{\weak}$ are in fact functions of the number
$\pi_0 = P\{\theta=0\}$, which we denote $\tilde g(\pi_0)$ and
$\tilde
g^{\weak}(\pi_0)$:
%
%
\begin{eqnarray}
\label{vf1} \tilde g(\pi_0) &=& \inf_{(h,\tau)} E
\bigl(C(h, \tau) \bigr),
\\
\tilde g^{\weak}(\pi_0) &=& \inf_{(h,\tau,X)} E
\bigl(C(h, \tau,X) \bigr). \label{vf1a}
\end{eqnarray}
Clearly, $\tilde g^{\weak} \leq\tilde g$. The following simple lemma
(see also \cite{PS2}, Section~2.7) provides important information about
the form of these two functions.

%
\begin{lemma}\label{lem0a}
The functions $\tilde g$ and $\tilde g^{\weak}$ are concave.
\end{lemma}

\begin{pf} By the law of total probability,
\begin{eqnarray*}
&&E \bigl(C(h, \tau) \bigr) \\
&&\qquad= \pi_0 E \biggl(a \tau(X) + b \int
^{\tau(X)}_0 h_t(X) \,dt \Big\vert\theta=0
\biggr)
\\
&&\qquad\quad{} + (1-\pi_0)\\
&&\qquad\qquad{}\times E \biggl(1_{\{\tau(X) < \theta\}} + a \bigl(\tau (X)-\theta
\bigr)1_{\{
\tau(X) > \theta\}} 
+ b \int^{\tau(X)}_0
h_t(X) \,dt \Big\vert\theta>0 \biggr).
\end{eqnarray*}
We note that the first expectation does not depend on $\pi_0$, since
$\tau(X)$ and $h_t(X)$ are determined by the observation process only,
and the second does not either, since the conditional distribution of
$\theta$ given that $\theta>0$ does not depend on $\pi_0$. Therefore,
$\pi_0 \mapsto E(C(h, \tau))$ is an affine function of $\pi_0$, and
$\tilde g$, being the infimum of affine functions, is concave. The same
argument applies to $\tilde g^{\weak}$.
\end{pf}

\subsection*{Sufficient statistic}

Let $\F^X = (\F^X_t)$ be the natural filtration of the observed
process $X$, augmented with $P$-null sets. Let $(\pi^h_t)$ be the
optional projection of $(1_{\{\theta\leq t \}}, t \geq0)$ onto this
filtration, so that for all $t$, $\pi^h_t = P\{\theta\leq t \mid X_s,
s \leq t\}$ a.s. The next several lemmas are identical both for
strategies and for control systems, so we state them only for strategies.

%
\begin{lemma}\label{lem0} With the above notation,
%
%
\begin{equation}
\label{cost0} E \bigl(C(h, \tau) \bigr) = E \biggl(1- \pi^h_\tau+
a \int^\tau_0 \pi^h_s
\,ds + b \int_0^\tau h_s \,ds
\biggr).
\end{equation}
\end{lemma}

\begin{pf} Note that $E(1_{\{\tau< \theta\}})= E(1-\pi^h_\tau)$ and
\begin{eqnarray*}
E \bigl((\tau-\theta)1_{\{\tau> \theta\}} \bigr)&=& E \biggl(\int
_0^\infty1_{\{\theta< s\}} 1_{\{s < \tau\}} \,ds
\biggr) = \int_0^\infty E \bigl(
\pi^h_s1_{\{s < \tau\}} \bigr) \,ds 
\\
&=& E
\biggl( \int_0^\tau\pi^h_s
\,ds \biggr).
\end{eqnarray*}
This proves the lemma.
\end{pf}

According to Lemma~\ref{lem0}, the expected cost associated to a
strategy $(h,\tau)$ is the expectation of an adapted functional of the
\emph{posterior probability} process $(\pi^h_t)$. Therefore, it will be
natural to express controls as functionals of $(\pi^h_t)$. We proceed
with the analysis of this process.

\section{Semimartingale characteristics of \texorpdfstring{$(\pi^h_t)$}{$(pi^h_t)$} and a
verification lemma}\label{sec3}

For $0\leq u < t$, let $\mu_{u,t}$ be the conditional distribution,
given that $\theta= u$, of $X$ restricted to $[0,t]$, and let $\mu_t$
be the unconditional distribution of $X$ restricted to $[0,t]$.

%
\begin{lemma}\label{rn} The Radon--Nikodym derivative of $\mu_{u,t}$
with respect
to $\mu_{t,t}$ is
%
%
\begin{equation}
\label{rnd} \frac{d \mu_{u,t}}{d \mu_{t,t}} = \exp \biggl(\int_u^t
\frac{r}{\sigma
^2} \,d X_s - \frac{1}{2}\int
_u^t \frac{r^2}{\sigma^2} h_s(X) \,ds
\biggr).
\end{equation}
\end{lemma}

\begin{pf} Recall Girsanov's theorem \cite{oksendal}, Theorem~8.6.6,
page 166: let
\begin{eqnarray*}
dZ_t &=& \sigma(Z_t) \,dW_t,
\\
d\tilde Z_t &=& \gamma_t \,dt + \sigma(\tilde
Z_t) \,dW_t,
\end{eqnarray*}
and suppose that under $P$, the process $(W_t)$ is a standard Brownian
motion. Define $\tilde P$ by
\[
\frac{d\tilde P}{dP} = \exp \biggl(- \int_0^t
\frac{\gamma_s}{\sigma
(\tilde Z_s)} \,dW_s - \frac{1}{2}\int_0^t
\biggl(\frac{\gamma_s}{\sigma(\tilde
Z_s)} \biggr)^2 \,ds \biggr).
\]
If $E_P(\frac{d\tilde P}{dP}) = 1$, then the law of $(\tilde Z_t)$
under $\tilde P$ is the same as the law of $(Z_t)$ under~$P$.

If $\theta= u$, then the law of $(X_s, s \leq t)$ is the same as
that of $(Y_s, s \leq t)$, where
%
%
\begin{equation}
\label{eYa} dY_s = r h_s(Y) 1_{\{u < s\}} \,ds +
\sigma\sqrt{h_s(Y)} \,dW_s,\qquad 0 < s < t.
\end{equation}
If $\theta= t$, then the law of $(X_s, s \leq t)$ is the same as that
of $(Z_s, s \leq t)$, where
\[
dZ_s = \sigma\sqrt{h_s(Z)} \,dW_s, \qquad 0 < s
< t.
\]
Therefore, for $A \in\B(C([0,t], \IR))$,
\[
\mu_{u,t} (A) = P\{Y_\cdot\in A\} = E_P
\bigl(1_A(Y_\cdot) \bigr) = E_{\tilde P}
\biggl(1_A(Y_\cdot) \frac{dP}{d\tilde P} \biggr),
\]
where $\tilde P$ is defined by
\begin{eqnarray*}
\frac{d\tilde P}{dP} &=& \exp \biggl(- \int_u^t
\frac{r h_s(Y)}{\sigma
\sqrt{h_s(Y)}} \,dW_s - \frac{1}{2}\int_u^t
\biggl(\frac{r h_s(Y)}{\sigma\sqrt
{h_s(Y)}} \biggr)^2 \,ds \biggr)
\\
&=& \exp \biggl(-\int_u^t \frac{r}{\sigma^2}
\sigma\sqrt{h_s(Y)} \,dW_s - \frac{1}{2}\int
_u^t \biggl(\frac{r}{\sigma}
\biggr)^2 h_s (Y) \,ds \biggr).
\end{eqnarray*}
Note in particular that Novikov's condition \cite{oksendal} is
satisfied. Using (\ref{eYa}), we see that this can be written
\[
\frac{d\tilde P}{dP} = \exp \biggl( - \int_u^t
\frac{r}{\sigma^2} \,dY_s + \frac{1}{2}\int
_u^t \biggl(\frac{r}{\sigma}
\biggr)^2 h_s(Y) \,ds \biggr).
\]

Therefore, by Girsanov's theorem,
\begin{eqnarray*}
\mu_{u,t}(A) &=& E_{\tilde P}\biggl(1_A(Z_\cdot)
\exp \biggl( \int_u^t \frac
{r}{\sigma^2}
\,dZ_s - \frac{1}{2}\int_u^t
\biggl( \frac{r}{\sigma} \biggr)^2 h_s(Z) \,ds \biggr)
\biggr)
\\
&=& \int_A \mu_{t,t} (d\omega) \exp \biggl( \int
_u^t \frac{r}{\sigma^2} \,dX_s -
\frac{1}{2}\int_u^t \biggl(
\frac{r}{\sigma} \biggr)^2 h_s(X) \,ds \biggr).
\end{eqnarray*}
This proves Lemma~\ref{rn}.
\end{pf}

Let $F_\theta$ denote the probability distribution function of $\theta
$, so that
\[
F_\theta(x) = \cases{ %
 0, &\quad
$\mbox{if } x < 0,$
\vspace*{2pt}\cr
\pi_0 + (1-\pi_0) \bigl(1-e^{-\lambda x} \bigr),&\quad
$\mbox{if } x \geq0.$}
\]

%
\begin{lemma}\label{lem3} We have
\[
\pi^h_t = \int_{0-}^t
\frac{d \mu_{u,t}}{d \mu_t} F_\theta(du) = \frac
{d\mu_{t,t}}{d\mu_t} \int
_{0-}^t \frac{d\mu_{u,t}}{d\mu_{t,t}} F_\theta(du)
\]
(note that the $0-$ accounts for the discontinuity of $F_\theta$ at $0$).
\end{lemma}

\begin{pf} The notation $\frac{d \mu_{u,t}}{d \mu_{t,t}}$ now refers to
the right-hand side of \eqref{rnd}, which is continuous in $u$. For the
first equality in the lemma, it suffices to show that for all $B \in\B
(C([0,t], \IR))$,
\[
E \biggl(1_{\{X\vert_{[0,t]} \in B \}} \int_{0-}^t
\frac{d \mu_{u,t}}{d
\mu_t} F_\theta(du) \biggr) = E (1_{\{X\vert_{[0,t]} \in B \}}
1_{\{
\theta\leq t \}} ).
\]
To see this, observe that
\begin{eqnarray*}
\int_{\{X\vert_{[0,t]}\in B\}} \,dP(\omega) \int_{0-}^t
\frac{d\mu
_{u,t}}{d\mu_t} (\omega) F_\theta(du) &=& \int_{0-}^t
F_\theta(du) \int_{\{X\vert_{[0,t]}\in B\}} \,dP(\omega)
\frac
{d\mu_{u,t}(\omega)}{d\mu_t}
\\
&=& \int_{0-}^t F_\theta(du)
\mu_{u,t} \{X\vert_{[0,t]} \in B \}
\\
&=& P \{ \theta\leq t, X\vert_{[0,t]} \in B \}.
\end{eqnarray*}
This proves the first equality. The second is a consequence of the
chain rule for Radon--Nikodym derivatives.
\end{pf}

%
\begin{lemma}\label{lem4} We have
\[
1 - \pi^h_t = \frac{d \mu_{t,t}}{d\mu_t} \int
_t^\infty F_\theta(du) = (1-
\pi_0) e^{-\lambda t} \frac{d \mu_{t,t}}{d\mu_t}.
\]
\end{lemma}

\begin{pf} As in Lemma~\ref{lem3}, one checks that
\[
P\{\theta> t \mid X_s, s \leq t\} = \int_t^{+\infty}
\frac{d \mu
_{u,t}}{d\mu_t} F_\theta(du).
\]
Since $\frac{d\mu_{u,t}}{d\mu_{t,t}} =1$ when $u>t$, the right-hand
side is equal to
\[
\int^{+\infty}_t \frac{d\mu_{u,t}}{d\mu_{t,t}}
\frac{d\mu_{t,t}}{d\mu
_t} F_\theta(du) = \frac{d\mu_{t,t}}{d\mu_t} \int
_t^{+\infty}F_\theta(du).
\]
This proves the first equality in the statement of the lemma. The
second equality is a consequence of the fact that for $u>0$, $F_\theta
(du) = (1-\pi_0)\lambda e^{-\lambda u} \,du$.
\end{pf}


Set
\[
\varphi^h_t = \frac{\pi^h_t}{1-\pi^h_t}
\]
and let
\[
Z_{u,t} = \int^t_u
\frac{r}{\sigma^2} \,dX_s - \frac{1}{2}\int
^t_u \frac
{r^2}{\sigma^2} h_s\, ds.
\]
Use Lemmas \ref{rn}, \ref{lem3} and \ref{lem4} to see that
%
%
\begin{eqnarray}
\label{ephi1}
\nonumber
\varphi^h_t &=&
\frac{e^{\lambda t}}{1-\pi_0} \int_{0-}^t \exp
(Z_{u,t} ) F_\theta(du)
\\
&=& \frac{e^{\lambda t}}{1-\pi_0} \exp(Z_{0,t} ) \int
_{0-}^t \exp(-Z_{0,u} )
F_\theta(du)
\\
&=& \frac{e^{\lambda t}}{1-\pi_0} \exp(Z_{0,t} ) \biggl(
\pi_0 + (1-\pi_0) \int_{0}^t
\exp(-Z_{0,u} )\lambda e^{-\lambda u} \,du \biggr).\nonumber
\end{eqnarray}

%
\begin{lemma}\label{lem5} The following s.d.e. is satisfied:
%
%
\begin{equation}
\label{ephi} d\varphi^h_t = \lambda \bigl(1+
\varphi^h_t \bigr) \,dt + \frac{r}{\sigma^2} \varphi
^h_t \,dX_t.
\end{equation}
\end{lemma}

\begin{pf} Observe from (\ref{e1.2}) that the quadratic variation of
$X_t$ is $d\langle X \rangle_t = \sigma^2 h_t \,dt$, so we can apply
It\^
o's formula and (\ref{ephi1}) to get
\begin{eqnarray*}
d\varphi^h_t &=& \lambda\varphi^h_t
\,dt
\\
&&{} + \frac{e^{\lambda t}}{1-\pi_0} \exp(Z_{0,t}) \biggl(
\frac{r}{\sigma
^2} \,dX_t - \frac{r^2h_t}{2 \sigma^2} \,dt + \frac{1}{2}
\biggl( \frac{r}{\sigma^2} \biggr)^2 \cdot\sigma^2
h_t \,dt \biggr)\\
&&\quad{}\times \int^t_{0-}
\exp(-Z_{0,u}) F_\theta(du)
\\
&&{} + e^{\lambda t} \exp(Z_{0,t})\exp(-Z_{0,t} ) \lambda
e^{-\lambda t} \,dt
\\
&=& \lambda \bigl(1+ \varphi^h_t \bigr) \,dt +
\frac{r}{\sigma^2} \varphi^h_t \,dX_t.
\end{eqnarray*}
\upqed\end{pf}

%
\begin{lemma}\label{lem6} The process $X=(X_t)_{t\ge0}$ has the
stochastic differential
\[
dX_t = r h_t(X) \pi_t \,dt + \sigma
\sqrt{h_t} \,d \bar W_t,
\]
where $(\bar W_t)$ is a standard Brownian motion.
\end{lemma}

\begin{pf} Observe that
\[
dX_t-rh_t \pi_t \,dt = (r h_t
1_{\{\theta\leq t\}} - r h_t \pi_t) \,dt + \sigma
\sqrt{h_t} \,dW_t,
\]
and the right-hand side has mean zero (given $X\vert_{[0,t]}$) and
quadratic variation $\sigma^2 h_t \,dt$. Further, the left-hand side is
adapted to $\F^X$, so that the right-hand side is too, and has mean
zero. In particular, it is the differential of a local $\F
^X$-martingale with quadratic variation $\sigma\sqrt{h_t} \,dt$.
According to \cite{KS}, Chapter~3, Theorem~4.2, this term is equal to
$\sigma\sqrt{h_t}$ times a standard Brownian motion increment. We note
for future reference that $(\bar W_t)$ need not be $\F^X$-adapted, but
the martingale
%
%
\begin{equation}
\label{emt} M_t = \sigma\int_0^t
\sqrt{h_s} \,d \bar W_t = X_t - \int
_0^t r h_s(X) \pi_s
\,ds
\end{equation}
is clearly $\F^X$-adapted.
\end{pf}
%

%
\begin{lemma}\label{lem7} Set $\rho= \frac{r}{\sigma}$. Then
%
%
\begin{equation}
\label{epiht1} d\pi^h_t = \lambda \bigl(1-
\pi^h_t \bigr) \,dt + \frac{r}{\sigma^2}
\pi^h_t \bigl(1-\pi^h_t \bigr)
\,dX_t - \frac{r^2}{\sigma^2} \bigl(\pi^h_t
\bigr)^2 \bigl(1-\pi^h_t \bigr)
h_t \,dt
\end{equation}
and
%
%
\begin{equation}
\label{epiht} d\pi^h_t = \lambda \bigl(1-
\pi^h_t \bigr) \,dt + \rho\pi^h_t
\bigl(1- \pi^h_t \bigr) \sqrt{h_t} \,d\bar
W_t.
\end{equation}
\end{lemma}

\begin{pf} Note that $\pi^h_t = \varphi^h_t (1+\varphi^h_t)^{-1} =
f(\varphi^h_t)$, where $f(x) = x(1+x)^{-1}$. Since $f^\prime(x) =
(1+x)^{-2}$ and $f^{\prime\prime}(x) = -2(1+x)^{-3}$, It\^o's formula
and Lemma~\ref{lem5} yield
\begin{eqnarray*}
d \pi^h_t &=& f^\prime \bigl(
\varphi^h_t \bigr) \,d \varphi^h_t
+ \frac{1}{2}f^{\prime
\prime} \bigl(\varphi^h_t
\bigr) \,d \bigl\langle\varphi^h \bigr\rangle_t
\\
&=& \frac{1}{(1+\varphi^h_t)^2} \biggl(\lambda \bigl(1+\varphi^h_t
\bigr) \,dt + \frac{r}{\sigma^2} \varphi^h_t
\,dX_t \biggr) + \frac{1}{2}\frac{-2}{(1+\varphi
^h_t)^3} \biggl(
\frac{r}{\sigma^2} \varphi^h_t \biggr)^2
\sigma^2 h_t \,dt.
\end{eqnarray*}
Recall that $1+\varphi^h_t = \frac{1}{1-\pi_t}$ to see that this is
equal to
\[
\lambda \bigl(1-\pi^h_t \bigr) \,dt + \frac{r}{\sigma^2}
\pi^h_t \bigl(1-\pi^h_t \bigr)
\,dX_t - \frac{r^2}{\sigma^2} \bigl(\pi^h_t
\bigr)^2 \bigl(1-\pi^h_t \bigr)
h_t \,dt,
\]
which establishes \eqref{epiht1}. By Lemma~\ref{lem6}, this is equal to
\[
\lambda \bigl(1-\pi^h_t \bigr) \,dt + \frac{r}{\sigma^2}
\pi^h_t \bigl(1-\pi^h_t \bigr)
\bigl(r h_t \pi^h_t \,dt + \sigma
\sqrt{h_t} \,d\bar W_t \bigr) - \frac{r^2}{\sigma
^2} \bigl(
\pi^h \bigr)^2_t \bigl(1-\pi^h_t
\bigr) h_t \,dt,
\]
which simplifies to
\[
\lambda \bigl(1- \pi^h_t \bigr) \,dt +
\frac{r}{\sigma} \pi^h_t \bigl(1-\pi^h_t
\bigr) \sqrt {h_t} \,d\bar W_t.
\]
This establishes \eqref{epiht}.
\end{pf}

\subsection*{Strategies expressed in terms of $(\pi^h_t)$}

According to \eqref{epiht}, $(\pi^h_t)$ is a diffusion process, and
therefore an optimal canonical control will typically be expressed as a
function of $\pi^h_t$; that is, we will mainly be interested in
controls $h_t(X)$ of the form $h_t(X) = h(t,\pi^h_t)$, where $h\colon
\IR_+ \times[0,1] \to[0,1]$ is measurable and given. We explain here
how to describe the observation process and the admissible control
$(h_t)$ associated with such a function $h$.

Consider the s.d.e.
%
%
\begin{eqnarray}
\label{e6.69ab0}
dp_t &=& \lambda(1-p_t) \,dt\nonumber \\
&&{}+
\frac{r}{\sigma^2} p_t (1-p_t) \bigl(r
h(t,p_t) 1_{\{\theta\leq t\}} \,dt + \sigma\sqrt{h(t,p_t)}
\,dW_t \bigr)
\\
&&{} - \frac{r^2}{\sigma^2}(p_t)^2 (1-p_t)
h(t,p_t) \,dt,\nonumber
\end{eqnarray}
with $p_0 = P\{\theta= 0\}$. Assume that $h$ is such that \eqref
{e6.69ab0} has a strong solution [i.e., an $(\F_t)$-adapted solution].
Then we define the observation process by $X_0 = 0$ and
%
%
\begin{equation}
\label{e6.69ab2} dX_t = r h(t,p_t) 1_{\{\theta\leq t\}} \,dt +
\sigma\sqrt{h(t,p_t)} \,dW_t.
\end{equation}
This process is adapted to $(\F_t)$, and by \eqref{e6.69ab0},
%
%
\begin{eqnarray}\label{e6.69ab1}
dp_t& =& \lambda(1-p_t) \,dt + \frac{r}{\sigma^2}
p_t (1-p_t) \,dX_t
\nonumber
\\[-8pt]
\\[-8pt]
\nonumber
&&{} - \frac
{r^2}{\sigma^2}(p_t)^2
(1-p_t) h(t,p_t) \,dt.
\end{eqnarray}
Let $q_t = p_t /(1-p_t)$. Applying It\^o's formula, we find that
%
%
\begin{equation}
\label{e6.69ab3} dq_t = \lambda(1+q_t) \,dt +
\frac{r}{\sigma^2} q_t \,dX_t.
\end{equation}
According to \cite{RY}, Chapter IX, (2.3), the solution of this linear
s.d.e. is
\begin{eqnarray*}
q_t &=& \exp \biggl(\frac{r}{\sigma^2} X_t + \lambda t -
\frac{1}{2}\frac
{r^2}{\sigma^4} \langle X \rangle_t 
 \biggr)\\ 
 &&{} \times
 \biggl[q_0 + \int_0^t
\exp \biggl(-\frac{r}{\sigma^2} X_s - \lambda s + \frac{1}{2}
\frac{r^2}{\sigma^4} \langle X \rangle_s 
 \biggr) \lambda
\,ds \biggr].
\end{eqnarray*}
In particular, $q_t$, and therefore $p_t$, is a function of $X\vert
_{[0,t]}$, and we can write $p_t = \hat h_t(X)$, where $(t,x) \mapsto
\hat h_t(x)$ from $\IR_+ \times C(\IR_+,\IR)$ to $\IR$ is progressively
measurable. Looking back to \eqref{e6.69ab2}, we see that $(X_t)$ is a
strong solution of the s.d.e.
%
%
\begin{equation}
\label{e6.69ab4} dX_t = r h_t(X) 1_{\{\theta\leq t\}} \,dt +
\sigma\sqrt{h_t(X)} \,dW_t,
\end{equation}
where $h_t(x) = h(t,\hat h_t(x))$. Therefore, $(h_t)$ is an admissible control.

Comparing \eqref{e6.69ab3} and \eqref{ephi}, we conclude that $q_t =
\varphi^h_t$ and therefore
%
%
\begin{equation}
\label{ept} p_t = \pi^h_t = P \bigl\{
\theta \leq t \mid\F^X_t \bigr\}.
\end{equation}
This means that the control $h_t(X)$ is indeed equal to $h(t,\pi^h_t)$.

We note that as in \eqref{epiht}, there is a Brownian motion $(\bar
W_t)$ such that
%
%
\begin{equation}
\label{e3.14} dp_t = \lambda(1-p_t) \,dt + \rho
p_t (1-p_t)\sqrt{h(t,p_t)} \,d\bar
W_t.
\end{equation}

If $\tau$ is a stopping time defined using $\pi^h_t$, for instance,
%
%
\begin{equation}
\label{tau0} \tau= \inf \bigl\{t \geq0\dvtx\pi^h_t \in
S \bigr\}
\end{equation}
for some Borel set $S \subset[0,1]$, then
\[
\tau= \inf \bigl\{t \geq0\dvtx\hat h_t(X) \in S \bigr\},
\]
so $\tau= \tau(X)$ is a canonical stopping time. 
In particular, $((h_t(x)),\tau(x))$ is a strategy.

The above discussion shows that if \eqref{e6.69ab0} has a strong
solution, then we can construct a strategy $((h_t),\tau)$ for which
\eqref{e1.2a0} or \eqref{e6.69ab4} admits a strong solution~$(X_t)$,
such that $p_t = \pi^h_t$, and the expected cost $E(C((h_t),\tau))$ is
given by \eqref{cost0}.

In the case where \eqref{e6.69ab0} admits a weak solution, we would
similarly conclude that \eqref{e1.2a} or \eqref{e6.69ab4} admits a weak
solution, and considering $\tau$ as in \eqref{tau0}, we would conclude
that $((h_t),\tau, X)$ is a control system with the same expected cost.

\subsection*{Verification lemma}

For $\pi\in[0, 1]$, let $E_\pi$ denote expectation in the case where
$\pi_0 = \pi$. 
Recall that we have defined
\[
\tilde g(\pi) = \inf_{(h, \tau)} E_\pi \bigl(C(h, \tau)
\bigr), \qquad\tilde g^{\weak}(\pi_0) = \inf_{(h,\tau,X)}
E \bigl(C(h, \tau,X) \bigr).
\]
By Lemma~\ref{lem0a}, $\tilde g$ is concave, and by Lemma~\ref{lem0},
\[
\tilde g(\pi) = \inf_{(h, \tau)} E_\pi \biggl(1-
\pi^h_\tau+ a \int_0^\tau
\pi^h_s \,ds + b \int_0^\tau
h_s \,ds \biggr),
\]
with similar properties for $\tilde g^{\weak}$.
According to \cite{elkaroui}, Theorem~3.67, we expect to be able to
characterize each of these two functions as a function $g^\ast$ with
certain properties concerning martingales and submartingales. The next
lemma gives conditions that will allow us to show that a function
$g^\ast$ is equal to $\tilde g$ (resp., $\tilde g^{\weak}$) and check
that a strategy $((h^\ast_t), \tau^\ast)$ [resp., a control system
$((h^\ast_t), \tau^\ast, X^\ast)$] is optimal.

%
\begin{lemma}[(Verification lemma)]\label{lem8}
Suppose that $g^\ast$ is a bounded continuous function defined on
$[0,1]$ such that $0 \leq g^\ast(x) \leq1 -x$, for all $x \in[0,1]$.
\begin{enumerate}[(1)]
\item[(1)] Suppose that for any $\pi\in[0,1]$, the following property holds:
\begin{enumerate}[(a)]
\item[(a)] for any strategy $((h_t), \tau)$ [resp., for any control
system $(h,\tau,X)$], the process $(Y_t)$ is an $\F^X$-submartingale
under $P_\pi$, where
%
%
\begin{equation}
\label{eY} Y_t = g^\ast \bigl(\pi_t^h
\bigr) + a \int_0^t \pi^h_s
\,ds + b \int_0^t h_s \,ds.
\end{equation}
\end{enumerate}
Then $g^\ast\leq\tilde g$ (resp., $g^\ast\leq\tilde g^{\weak}$).

\item[(2)] Suppose that for any $\pi\in[0,1]$, in addition to \textup{(a)}, the
following three properties hold:
\begin{enumerate}[(a)]
\item[(b)] for the strategy $((h^\ast_t), \tau^\ast)$ [resp., the
control system $((h^\ast_t), \tau^\ast, X^\ast)$], the process
$(Y^\ast
_{t \wedge\tau^\ast})$ is an $\F^X$-martingale under $P_\pi$, where
\[
Y^\ast_t = g^\ast \bigl(\pi^{h^\ast}_t
\bigr) + a \int_0^t \pi_s^{h^\ast}
\,ds + b \int_0^t h^\ast_s
\,ds;
\]

\item[(c)] $E_\pi(\tau^\ast) < + \infty$;

\item[(d)] $g^\ast(\pi^{h^\ast}_{\tau^\ast}) = 1- \pi^{h^\ast
}_{\tau
^\ast}$.
\end{enumerate}
Then $g^\ast= \tilde g$ and $((h^\ast_t), \tau^\ast)$ is an
optimal strategy [resp., $g^\ast= \tilde g^{\weak}$, and $((h^\ast_t),
\tau^\ast, X^\ast)$ is an optimal control system].
\end{enumerate}
\end{lemma}

\begin{pf} We first establish (1). Let $((h_t), \tau)$ be a strategy.
If $E(\tau) = + \infty$, then $E(C(h, \tau)) = + \infty$. Indeed, by
(\ref{cost}), $E(C(h, \tau)) \geq a E(\tau1_{\{\tau> \theta\}}) -
a E(\theta)$. Since
\[
E(\tau) = E(\tau1_{\{\tau> \theta\}}) + E(\tau1_{\{\tau\leq
\theta
\}})
\]
and the second term is no greater than $E(\theta) < + \infty$, we
conclude that $E(\tau1_{\{\tau> \theta\}}) = + \infty$ and so
$E(C(h, \tau)) = + \infty$.

Therefore, in the definition of $\tilde g$, we can restrict the
infimum to those strategies for which $E(\tau) < + \infty$. Since $1-x
\geq g^\ast(x)$, Lemma~\ref{lem0} implies that
\[
E_\pi \bigl(C(h, \tau) \bigr) \geq E(Y_\tau).
\]
Since $(Y_t)$ is a submartingale by (a) and $t \wedge\tau$ is a
bounded stopping time, $E_\pi(Y_{t \wedge\tau}) \geq E_\pi(Y_0) =
g^\ast(\pi)$. By Fatou's lemma in the form $E(\limsup Y_{n\wedge\tau})
\geq\limsup E(Y_{n\wedge\tau})$ (cf. \cite{CD}, Chapter~1), which
applies since $E(\tau) < + \infty$, we see that
\[
E_\pi(Y_\tau) \geq\limsup_{t \to\infty}
E_\pi(Y_{t \wedge\tau}) \geq g^\ast(\pi).
\]
We conclude that $E_\pi(C(h, \tau)) \geq g^\ast(\pi)$ for all
strategies $((h_t), \tau)$, and therefore $\tilde g \geq g^\ast$. The
proof for $\tilde g^{\weak}$ is identical and is omitted.

We now establish (2) for $\tilde g$. It suffices to show that $g^\ast
(\pi) = E_\pi(Y_{\tau^\ast}^\ast)$. Indeed, this will complete the
proof, since by (d) and Lemma~\ref{lem0},
\begin{eqnarray*}
g^\ast(\pi) &=& E_\pi \bigl(Y^\ast_{\tau^\ast}
\bigr) = E_\pi \biggl(g^\ast \bigl(\pi^{h^\ast}_{\tau^\ast}
\bigr) + a \int_0^{\tau^\ast} \pi^{h^\ast}_s
\,ds + b \int_0^{\tau^\ast} h^\ast_s
\,ds \biggr)
\\
&=& E_\pi \bigl(C \bigl(h^\ast, \tau^\ast \bigr)
\bigr) \geq\tilde g(\pi).
\end{eqnarray*}
Since we have already proved that $\tilde g \geq g^\ast$, this shows
that $g^\ast(\pi) = \tilde g(\pi)$.

In order to check that $g^\ast(\pi) = E_\pi(Y^\ast_{\tau^\ast})$, note
that $0 \leq Y^\ast_t \leq1 + (a+b)t$ and $E_\pi(\tau^\ast) < +
\infty
$ by (c). Therefore, $(Y^\ast_{t \wedge\tau^\ast})$, which is a
martingale by (b), is uniformly integrable. By the optional sampling
theorem \cite{durrett}, $E(Y_\tau^\ast) = E(Y^\ast_0) = g^\ast(\pi)$.
This completes the proof for $\tilde g$. The proof for $\tilde g^{\weak
}$ is identical and is omitted.
\end{pf}
%

\section{A candidate for the value function}\label{sec4}

We now seek analytical conditions on a function $g^\ast$ that will
guarantee the properties of Lemma~\ref{lem8}. Consider the process
$(Y_t)$ defined in (\ref{eY}) (we write $g$ instead of $g^\ast$ to
simplify the notation). By It\^o's formula and Lemma~\ref{lem7},
%
%
\begin{eqnarray}\label{drift}
\nonumber
dY_t &=& g^\prime \bigl(\pi^h_t
\bigr) \,d\pi^h_t + \frac{1}{2}g^{\prime\prime}
\bigl( \pi^h_t \bigr) \,d \bigl\langle\pi^h
\bigr\rangle_t + a \pi^h_t \,dt + b
h_t \,dt
\\
 &=& \biggl[\lambda g^\prime \bigl(\pi^h_t
\bigr) \bigl(1-\pi^h_t \bigr) + \frac{1}{2}g^{\prime
\prime}
\bigl(\pi^h_t \bigr) \bigl(\rho\pi^h_t
\bigl(1-\pi^h_t \bigr) \bigr)^2
h_t + a \pi^h_t + bh_t \biggr]
\,dt
\\
\nonumber
&&{} + g^\prime \bigl(\pi^h_t \bigr)
\frac{r}{\sigma} \pi^h_t \bigl(1-\pi^h_t
\bigr) \sqrt{h_t} \,d\bar W_t.
\end{eqnarray}
Therefore, $(Y_t)$ will be a submartingale if the term in brackets is
nonnegative, for any value of $h_t$. Since this term is an affine
function of $h_t$, this is equivalent to this term being nonnegative
for $h_t = 0$ and $h_t = 1$, that is, for all $x \in[0,1]$,
%
%
\begin{equation}
\label{e1.4a} \lambda g^\prime(x) (1-x) + a x \geq0
\end{equation}
and
%
%
\begin{equation}
\label{e1.4b} \lambda g^\prime(x) (1-x) + \tfrac{1}{2}g^{\prime
\prime}(x)
\bigl(\rho x(1-x) \bigr)^2 + a x + b \geq0.
\end{equation}

\subsection*{Intuition and smooth fit}

We can imagine that the optimal strategy, in either the strong or the
weak formulation, is of the following form: do not observe if $\pi^h_t$
is small, declare the alarm if $\pi^h_t$ is close to 1 and observe
otherwise. More precisely, we postulate that there are two constants $0
\leq A \leq B \leq1$ such that on $[0, A]$, it is optimal \emph{not}
to observe, on $]A, B[$ it is optimal to observe without declaring an
alarm and on $[B,1]$, it is optimal to stop and declare the alarm. That is,
%
%
\begin{equation}
\label{e1.4c} h^\ast_t = 1_{\{\pi^{h^\ast}_t > A\}} \quad\mbox{and}\quad
\tau^\ast= \inf \bigl\{t \geq0\dvtx\pi^{h^\ast}_t \geq
B \bigr\}.
\end{equation}
In order to satisfy condition (b) of Lemma~\ref{lem8}, we need
%
%
\begin{equation}
\label{e1.5} \lambda g^\prime(x) (1-x) + ax = 0,\qquad x \in\,]0,A]
\end{equation}
and
%
%
\begin{equation}\qquad
\label{e1.6} \lambda g^\prime(x) (1-x) + ax + \tfrac{1}{2}g^{\prime
\prime}(x)
\rho^2 x^2(1-x)^2 + b= 0,\qquad x \in\,]A,B[.
\end{equation}
In order to satisfy condition (d) of Lemma~\ref{lem8}, we need
%
%
\begin{equation}
\label{e1.7} g(x) = 1-x,\qquad  x \in[B,1].
\end{equation}

In order to find an expression for $g$, it is natural to solve first
the differential equations (\ref{e1.5}) and (\ref{e1.6}) separately,
that is, to seek two functions $g_1$ and $g_2$ such that
%
%
\begin{equation}
\label{e1.8} \lambda g^\prime_1(x) (1-x) + ax = 0,\qquad 0 < x
< A
\end{equation}
and
%
%
\begin{equation}\qquad
\label{e1.9} \lambda g^\prime_2(x) (1-x) + ax +
\tfrac{1}{2} g^{\prime\prime}_2(x) \rho^2
x^2(1-x)^2 + b = 0,\qquad  A < x < B.
\end{equation}
Three constants of integration will appear, one for $g_1$ and two for
$g_2$. These constants can then be determined by ``pasting together''
$g_1$ and $g_2$, that is, requiring equalities such as
%
%
\begin{equation}
\label{e1.10} g_1(A) = g_2(A)
\end{equation}
and, by (\ref{e1.7}),
%
%
\begin{equation}
\label{e1.11} g_2(B) = 1-B.
\end{equation}
These two equalities are referred to as ``continuous fit'' \cite{PS}. As
in most problems of optimal stopping or control, they are not
sufficient to determine the five unknown constants, namely, the three
constants of integration and the two ``free boundaries'' $A$ and $B$.
For this, it is necessary to use a version of the ``principle of smooth
fit''; see \cite{PS}. In particular, one can postulate that
%
%
\begin{equation}
\label{e1.12} g^\prime_2(B) = -1
\end{equation}
and
%
%
\begin{equation}
\label{e1.13} g^\prime_1(A) = g^\prime_2(A).
\end{equation}
We need one more equation in addition to (\ref{e1.10})--(\ref{e1.13}),
since there are five unknown constants. Since we want to apply It\^o's
formula, it is natural to want $g$ to be twice differentiable at $A$.
This gives one more equation,
%
%
\begin{equation}
\label{e1.13a} g^{\prime\prime}_1(A) = g^{\prime\prime}_2(A).
\end{equation}



\subsection*{Solving the equations}

We seek functions $g_1$ and $g_2$ defined on $[0,1]$ satisfying (\ref
{e1.8})--(\ref{e1.13a}). Set
\[
f_1 (x) = g^\prime_1(x),\qquad f_2 (x) =
g^\prime_2(x).
\]

\subsection*{The value of $A$}

For $0<x<A$, differentiate (\ref{e1.8}) to get
\[
-\lambda f_1(x) + \lambda f^\prime_1(x) (1-x)
+ a = 0,
\]
that is,
%
%
\begin{equation}
\label{e1.14} f^\prime_1(x) = \frac{\lambda f_1(x)-a}{\lambda(1-x)}.
\end{equation}
From (\ref{e1.9}), we get
%
%
\begin{equation}
\label{e1.15} f^\prime_2(x) = \frac{-ax-b-\lambda f_2(x)(1-x)}{({1}/{2})
\rho^2 x^2 (1-x)^2}.
\end{equation}
By (\ref{e1.13a}), if we plug $x = A$ into (\ref{e1.14}), (\ref
{e1.15}), we get
\[
- \lambda(aA+b) - \lambda^2 f_2(A) (1-A) = -
\frac{a}{2} \rho^2 A^2(1-A) + \frac{\lambda\rho^2}{2}
A^2(1-A) f_1(A).
\]
Since $f_2(A) = f_1(A)$ by (\ref{e1.13}), we solve for $f_1(A)$,
%
%
\begin{equation}
\label{e1.16} f_1(A) = \frac{({a\rho^2}/{2}) A^2(1-A) - \lambda
(aA +
b)}{(1-A)(\lambda^2 + ({\lambda\rho^2}/{2}) A^2)}.
\end{equation}
Plugging (\ref{e1.16}) into (\ref{e1.8}) gives an equation for $A$,
whose solution is
%
%
\begin{equation}
\label{e1.17} A = \sqrt{\frac{2 \lambda b}{a \rho^2}}.
\end{equation}

For the observation region $]A,B[$ to be nonempty, we must have $A <
1$, but further, since we want $g_1$ to be concave by Lemma~\ref
{lem0a}, we also must have
%
%
\begin{equation}
\label{e1.18} f_1(A) = g^\prime_1(A) > -1.
\end{equation}
From (\ref{e1.8}),
%
%
\begin{equation}
\label{e1.19} g^\prime_1(x) = - \frac{a}{\lambda}
\frac{x}{1-x},
\end{equation}
so (\ref{e1.8}) and (\ref{e1.18}) give
$
- \frac{a}{\lambda} \frac{A}{1-A} > -1$,
or equivalently,
$
A < \frac{\lambda}{a+\lambda}$.
With (\ref{e1.17}), we conclude that the observation region $]A,B[$ is
not empty if
%
%
\begin{equation}
\label{e1.20} b < \frac{\lambda a \rho^2}{2(a+\lambda)^2}.
\end{equation}

\subsection*{Determining $f_2(x)$}

For $A<x<B$, equation (\ref{e1.9}) becomes
%
%
\begin{equation}
\label{e1.21} \lambda f_2(x) (1-x) + ax + \tfrac{1}{2}f^\prime_2(x)
\rho^2 x^2(1-x)^2 + b = 0.
\end{equation}
A solution of the homogeneous equation
\[
\lambda f(x) (1-x) + \tfrac{1}{2}f^\prime(x) \rho^2
x^2(1-x)^2 = 0
\]
is
%
%
\begin{equation}
\label{e1.21a} f(x) = \biggl(\frac{1-x}{x} \biggr)^\alpha
e^{\alpha/x} \qquad\mbox{where } \alpha= \frac{2 \lambda}{\rho^2}.
\end{equation}
Therefore, the solution of the inhomogeneous equation (\ref{e1.21}) is
%
%
\begin{equation}
\label{e1.22} f_2(x) = K_1 f(x) + f(x) \int
_A^x \frac{-2}{\rho^2} \frac{ay +
b}{y^2(1-y)^2}
\frac{1}{f(y)} \,dy.
\end{equation}
From (\ref{e1.13}) and (\ref{e1.8}), we conclude that
%
%
\begin{equation}
\label{e1.23} K_1 = - \frac{a}{\lambda} \frac{A}{1-A}
\frac{1}{f(A)}.
\end{equation}
Formulas (\ref{e1.23}) and (\ref{e1.22}) together determine $f_2(x)$.

%
\begin{remark}{ In the case where $b=0$, then $A=0$ by \eqref{e1.17},
and we must have $K_1=0$ in order that $f_2(x)$ be bounded. This
recovers the case discussed in~\cite{shiryaev}, Chapter~4.4. Therefore,
we consider the case $b>0$.
}
\end{remark}

\subsection*{Determining B}

Observe that
\[
\lim_{x \to1} f(x) = 0\quad \mbox{and}\quad \lim_{x \to1}
f_2(x) = -\infty.
\]
Indeed, the first equality is obvious, and the second holds because for
$x$ near~1,
\[
f(x) \sim(1-x)^\alpha,
\]
and, using l'Hopital's rule,
\begin{eqnarray*}
f_2(x) &\sim& \frac{-2}{\rho^2} (1-x)^\alpha\int
_A^x \frac
{(a+b)e^{-\alpha}}{(1-y)^{2+\alpha}} \,dy
\sim-(1-x)^{1+\alpha} (1-x)^{-2-\alpha}
\\
&\sim& -(1-x)^{-1}.
\end{eqnarray*}
Therefore, if (\ref{e1.20}) holds, then $f_2(A) = K_1 f(A) = - \frac
{a}{\lambda} \frac{A}{1-A} > -1$, so there is $B \in\,]A, 1[$ such that
%
%
\begin{equation}
\label{e1.24} f_2(B) = -1.
\end{equation}
%
With this choice of $B$, (\ref{e1.12}) is satisfied. The next lemma
shows that in fact, there is only one solution to \eqref{e1.24}.

%
\begin{lemma} \label{lem4.2}
The function $f_2$ defined in \eqref{e1.22} is strictly decreasing on
$[A,1[$, and therefore, there is a unique $B \in\,]A,1[$ satisfying
\eqref{e1.24}.
\end{lemma}

\begin{pf} By \eqref{e1.21},
%
%
\begin{equation}
f_2'(x) = \frac{2\lambda}{\rho^2} \frac{1}{x^2(1-x)} \bigl(
\psi(x) - f_2(x) \bigr),
\end{equation}
where
\[
\psi(x) = -\frac{ax+b}{\lambda(1-x)}.
\]
Therefore, $f_2'(x) < 0$ if and only if $\psi(x) < f_2(x)$. In fact, we
will see in \eqref{e1.35aa} [see also \eqref{e1.27a}] that
\[
f_2(x) > -\frac{a}{\lambda} \frac{x}{1-x} > \psi(x),\qquad x \in\,]A,1[.
\]
We conclude that $f_2'(x) <0$ for $x \in\,]A,1[$, and this proves the lemma.
\end{pf}
%

%
\begin{remark} It is clear from \eqref{e1.17} that the value of the
constant $A$ is a continuous function of the observation cost $b$. The
same is true for the constant $B$, by the following argument.

We make explicit the dependence of $f_2$ on $b$ by writing $f_2(x,b)$.
Equation~\eqref{e1.24} becomes $f_2(B,b) = -1$. We see that $\frac
{\partial f_2}{\partial b}(x,b) > 0$ by differentiating under the
integral sign in \eqref{e1.22}. Therefore, the implicit function
theorem implies that $B$ is a continuous (and even differentiable)
function of~$b$.
\end{remark}

\subsection*{Determining $g_2(x)$}

Because $g^\prime_2(x) = f_2(x)$, $g_2(x)$ can be written
%
%
\begin{equation}
\label{e1.25} g_2(x) = \int_A^x
f_2(y) \,dy + K_2.
\end{equation}
From (\ref{e1.11}), we see that
%
%
\begin{equation}
\label{e1.26} K_2 = 1-B - \int_A^B
f_2(y) \,dy,
\end{equation}
so that
%
%
\begin{equation}
\label{e1.26a} g_2(x) = \int_B^x
f_2(y) \,dy + 1-B.
\end{equation}

\subsection*{Determining $g_1(x)$}

Because $g^\prime_1(x) = f_1(x)$, $g_1(x)$ can be written
%
%
\begin{equation}
\label{e1.27} g_1(x) = \int_A^x
f_1(y) \,dy + K_3,
\end{equation}
where $f_1(x)$ is determined from (\ref{e1.8}),
%
%
\begin{equation}
\label{e1.27a} f_1(x) = - \frac{a}{\lambda} \frac{x}{1-x}.
\end{equation}
From (\ref{e1.10}), (\ref{e1.25}) and (\ref{e1.27}), we get
%
%
\begin{equation}
\label{e1.28} K_3 = K_2.
\end{equation}
We can perform the integration in (\ref{e1.27}) to get
%
%
\begin{equation}
\label{e1.28a} g_1(x) = \frac{a}{\lambda} \bigl(x + \ln(1-x) - A -
\ln(1-A) \bigr) + K_2,
\end{equation}
with $K_2$ determined by (\ref{e1.26}).

We have now found two functions $g_1$ and $g_2$ that solve \eqref
{e1.8}--\eqref{e1.13a}. In order to ensure that this solves our optimal
control problem, slightly more is needed: in particular, we need
inequalities (\ref{e1.4a}) and (\ref{e1.4b}) for all $x\in[0,1]$. Set
%
%
\begin{eqnarray}
\label{e1.28b} L_1 g(x) &=& \lambda g^\prime(x) (1-x) + ax,
\\
L_2 g(x) &=& \lambda g^\prime(x) (1-x) + ax +
\tfrac{1}{2} g^{\prime\prime
}(x) \rho^2 x^2(1-x)^2
+ b. \label{e1.28c}
\end{eqnarray}


%
\begin{prop}[(Candidate value function)]\label{lem10} Suppose that $0<
b < \lambda a
\rho^2 / (2(a+\lambda)^2)$. Define $g(x)$ on $[0,1]$ by
%
%
\begin{equation}
\label{defg} g(x) = \cases{ %
g_1(x), &\quad $\mbox{if } 0 \leq x \leq A,$
\vspace*{2pt}\cr
g_2(x), &\quad $\mbox{if } A \leq x \leq B,$
\vspace*{2pt}\cr
1-x, &\quad $\mbox{if } B \leq x \leq1,$}
\end{equation}
where $A$ is defined in \eqref{e1.17}, and $B$ is defined in \eqref
{e1.24}. Then $g$ is strictly concave in $[0,B]$, and
%
%
\begin{eqnarray}
\label{e1.29a} 0 &\leq& g(x) \leq1-x, \qquad 0 \leq x \leq1,
\\
\label{e1.29aa} L_1 g(x) &=& 0,\qquad 0 \leq x \leq A,
\\
\label{e1.29aaa} L_2 g(x) &=& 0, \qquad A \leq x < B.
\end{eqnarray}
Furthermore,
%
%
\begin{eqnarray}
\label{e1.29} L_2g(x) &\geq&0, \qquad 0 \leq x \leq A,
\\
\label{e1.30} L_1 g(x) &\geq&0, \qquad A \leq x \leq B,
\\
\label{e1.31} L_1g(x) &\geq&0, \qquad B \leq x \leq1,
\\
\label{e1.32} L_2g(x) &\geq&0, \qquad B \leq x \leq1.
\end{eqnarray}
\end{prop}

\begin{pf} Properties \eqref{e1.29aa} and \eqref{e1.29aaa} follow from
the construction of $g_1$ and~$g_2$; see \eqref{e1.8} and \eqref{e1.9}.
The strict concavity of $g_1$ and $g_2$ (hence of $g$ on $[0,B]$)
follow from \eqref{e1.27a} and Lemma~\ref{lem4.2}. This concavity
property and \eqref{e1.24} imply $g(x) \leq1-x$, $0 \leq x \leq1$.
Finally, since $g_1'(x) = f_1(x) <0$ for $0 < x \leq A$ and $g_2'(x) =
f_2(x) \leq0$ for $A \leq x \leq B$, $g$ is nondecreasing on $[0,B]$,
therefore nonnegative on $[0,B]$ since $g_2(B) = 1-B \geq0$. This
proves \eqref{e1.29a}.

Note that (\ref{e1.31}) implies (\ref{e1.32}), and on $[B,1]$, (\ref
{e1.31}) becomes
$
- \lambda(1-x) + ax \geq0$,
that is,
$
x \geq\frac{\lambda}{a+\lambda}$.
Therefore, (\ref{e1.31}) will hold provided we show that
%
%
\begin{equation}
\label{e1.32a} B \geq\frac{\lambda}{a+\lambda}.
\end{equation}

To see this, note from (\ref{e1.27a}) that
%
%
\begin{equation}
\label{e1.32b} f_1 \biggl(\frac{\lambda}{a+\lambda} \biggr) = -
\frac{a}{\lambda}\cdot\frac{{\lambda/}{(a+ \lambda)}}{1- {\lambda}/{(a+\lambda)}} = -1.
\end{equation}
We shall show that
%
%
\begin{equation}
\label{e1.33} f_2(x) \geq f_1(x) \qquad\mbox{for } x \geq A.
\end{equation}
Then, (\ref{e1.33}) and (\ref{e1.32b}) imply that
\[
f_2 \biggl(\frac{\lambda}{a+\lambda} \biggr) \geq-1\qquad \mbox{that is, } B
\geq\frac{\lambda}{a+\lambda}
\]
[since $f_2(x) < -1$ for $x > B$, by \eqref{e1.24} and Lemma~\ref
{lem4.2}], proving (\ref{e1.32a}).

It remains to prove (\ref{e1.33}). Set $h(x) = f_2(x) - f_1(x)$. From
(\ref{e1.8}) and (\ref{e1.9}), we see that for $x > A$,
%
%
\begin{equation}\qquad
\label{e1.34} \lambda h(x) (1-x) + \tfrac{1}{2}h^\prime(x)
\rho^2 x^2(1-x)^2 + b + \tfrac{1}{2}f^\prime_1(x)
\rho^2 x^2(1-x)^2 = 0.
\end{equation}
By (\ref{e1.19}),
%
%
\begin{equation}
\label{e1.34a} f^\prime_1(x) = - \frac{a}{\lambda}
\frac{1}{(1-x)^2},
\end{equation}
so (\ref{e1.34}) becomes
%
%
\begin{equation}
\label{e1.35} \lambda h(x) (1-x) + \frac{1}{2}h^\prime(x)
\rho^2 x^2(1-x)^2 + b - \frac
{a \rho^2}{2 \lambda}
x^2 = 0.
\end{equation}

Recall from (\ref{e1.17}) that $b - \frac{a \rho^2}{2 \lambda} x^2 <
0$ for $x > A$. We note that $h(A) = h^\prime(A) = 0$ by (\ref{e1.13})
and (\ref{e1.13a}), and from (\ref{e1.35}), the following holds:
for $x > A$, it is not possible to have simultaneously $h(x) < 0$ and
$h^\prime(x) < 0$. Since $h(A) = 0$, this implies that for $x > A$,
$h(x)$ cannot be negative (since otherwise, there would be $y \in\,
]A,x[$ with $h(y) < 0$ and $h^\prime(y) <0$), therefore $h(x) > 0$ for
$x > A$, that is,
%
%
\begin{equation}
\label{e1.35aa} f_2(x) > f_1(x) \qquad\mbox{for } x > A.
\end{equation}
This proves (\ref{e1.33}). Therefore, (\ref{e1.31}) is proved.

To check (\ref{e1.29}), we use (\ref{e1.8}), to see that for $0 \leq x
\leq A$,
\[
L_2g(x) = \tfrac{1}{2}g^{\prime\prime}_1(x)
\rho^2 x^2(1-x)^2 + b,
\]
and from (\ref{e1.34a}),
\[
g_1^{\prime\prime}(x) = - \frac{a}{\lambda} \frac{1}{(1-x)^2},
\]
therefore,
\[
L_2g(x) = - \frac{a}{2 \lambda} \rho^2 x^2 +
b, \qquad x \leq A,
\]
and the right-hand side is nonnegative for $x \leq A$ by (\ref{e1.17}).
This proves (\ref{e1.29}).

Finally, (\ref{e1.30}) is a consequence of (\ref{e1.33}), since (\ref
{e1.33}) implies that
\[
L_1g_2(x) \geq L_1 g_1(x) =
0.
\]
\upqed\end{pf}

\subsection*{Case where $b \geq\frac{\lambda a \rho^2}{2(a+\lambda)^2}$}

In this case, we postulate that the observation region $]A, B[$ is
empty (i.e., \mbox{$B =A)$}, so we seek $g_1(x)$ such that
%
%
\begin{eqnarray}
\label{e1.35a} \lambda g^\prime_1(x) (1-x) + a x &=& 0,\qquad 0
\leq x < B,\\
\label{e1.35b} g_1 (B) &=& 1-B,
\\
\label{e1.35c} g^\prime_1(B) &=& -1.
\end{eqnarray}
From (\ref{e1.35a}), we see that
%
%
\begin{equation}
\label{e1.36} g_1^\prime(x) = - \frac{a}{\lambda}
\frac{x}{1-x} = \frac{a}{\lambda} \biggl(1- \frac{1}{1-x} \biggr),
\end{equation}
so for some constant $K$ to be determined,
%
%
\begin{equation}
\label{e1.37} g_1(x) = K + \frac{a}{\lambda} x +
\frac{a}{\lambda} \ln(1-x).
\end{equation}
From (\ref{e1.36}) and (\ref{e1.35c}), we see that
\[
\frac{a}{\lambda} \biggl(1- \frac{1}{1-B} \biggr) = -1,
\]
that is,
%
%
\begin{equation}
\label{e1.37a0} B = \frac{\lambda}{a+\lambda}.
\end{equation}
From (\ref{e1.35b}) and (\ref{e1.37}), we obtain
\[
K = 1- B - \frac{a}{\lambda} B - \frac{a}{\lambda} \ln(1-B)
= -
\frac{a}{\lambda} \ln \biggl(\frac{a}{a+\lambda} \biggr).
\]
Therefore,
%
%
\begin{equation}
\label{e1.37a} g_1(x) = \frac{a}{\lambda} x + \frac{a}{\lambda}
\biggl(\ln(1-x) - \ln \biggl(\frac{a}{a+\lambda} \biggr) \biggr).
\end{equation}
We note that $g_1^\prime(x)$ is decreasing, $g_1^\prime(0) = 0$ and
$g^\prime_1(B) = -1$, so $1-x \geq g_1(x)$ for $0 \leq x \leq B$, by
(\ref{e1.35b}). Since $1-x \geq a/(a + \lambda) = 1-B$ for $x \leq B$,
$g_1(x) \geq0$ for $0\leq x \leq B$.

%
\begin{prop}[(Candidate value function)]\label{lem11}  Suppose that $b
\geq\lambda a
\rho^2 /\break (2(a+\lambda)^2)$. Define $g_1(x)$ as in \eqref{e1.37a} and
$g(x)$ on $[0,1]$ by
%
%
\begin{equation}
\label{defg2} g(x) = \cases{ %
g_1(x), & \quad$\mbox{if } 0 \leq x \leq B,$
\vspace*{2pt}\cr
1-x, &\quad $\mbox{if } B \leq x \leq1,$}
\end{equation}
where $B$ is defined in \eqref{e1.37a0}. Then $g$ is strictly concave
on $[0,B]$,
%
%
\begin{eqnarray}
\label{e1.38a} 0 &\leq&g(x) \leq1-x, \qquad 0 \leq x \leq1,
\\
\label{e1.38aa} L_1 g(x) &=&0, \qquad x \in[0,B],
\end{eqnarray}
and furthermore,
%
%
\begin{eqnarray}
\label{e1.38} L_1 g(x) &\geq&0\qquad \mbox{for } B \leq x \leq1,
\\
\label{e1.39} L_2 g(x) &\geq&0\qquad\mbox{for } 0 \leq x \leq1.
\end{eqnarray}
\end{prop}

\begin{pf} Property \eqref{e1.38aa} follows from \eqref{e1.35a}, and
the strict concavity of $g_1$, hence of $g$, on $[0,B]$ and \eqref
{e1.38a} are established just after \eqref{e1.37a}.

Note that for $B \leq x \leq1$,
\[
L_1 g(x) \geq0 \quad\Longleftrightarrow\quad- \lambda(1-x) + a x \geq0
\quad\Longleftrightarrow\quad x \geq\frac{\lambda}{a+\lambda} = B,
\]
and this is indeed that case, so (\ref{e1.38}) holds.

For $B \leq x \leq1$, $L_2 g(x) = L_1g(x) + b$, and both of these
terms are nonnegative, so $L_2 g(x) \geq0$ for these $x$, proving part
of (\ref{e1.39}).

For $0 < x < B$,
\[
L_2 g(x) = L_2 g_1(x) = L_1
g_1(x) + \tfrac{1}{2}g_1^{\prime\prime}(x) \rho
^2 x^2(1-x)^2 + b.
\]
Since $L_1 g_1(x) = 0$,
\[
L_2 g(x) = \frac{1}{2}\frac{a}{\lambda} \frac{-1 }{(1-x)^2}
\rho^2 x^2 (1-x)^2 + b = - \frac{a \rho^2 x^2}{2 \lambda}
+ b,
\]
so
\[
L_2g(x) \geq0 \quad\Longleftrightarrow\quad\frac{a \rho^2 x^2}{2 \lambda} \leq b
\quad\Longleftrightarrow\quad x \leq\sqrt{\frac{2 \lambda b}{a \rho^2}}.
\]
This will hold for $x \leq B$ provided it holds for $x=B$. Now
\[
B \leq\sqrt{\frac{2 \lambda b}{a \rho^2}}\quad \Longleftrightarrow\quad \biggl(\frac{\lambda}{a+\lambda}
\biggr)^2 \leq\frac{2 \lambda b}{a \rho^2} \quad\Longleftrightarrow\quad b \geq
\frac{\lambda a \rho^2}{2(a+\lambda)^2},
\]
which is the assumption of this case. This proves (\ref{e1.39}).
\end{pf}

\subsection*{Comments on the optimal strategy}

In the case where $b \geq\lambda a \rho^2 / (2(a + \lambda)^2)$, the
observation region is empty, the candidate optimal control is $h_t^*
\equiv0$ [with this control, \eqref{e1.2a} obviously has a strong
solution] and the candidate optimal stopping time is
%
%
\begin{equation}
\label{en1a} \tau^* = \inf \bigl\{ t \geq0\dvtx\pi_t^* \geq B
\bigr\},
\end{equation}
where $(\pi_t^*)$ is defined by
%
%
\begin{equation}
\label{en1} d \pi_t^* = \lambda \bigl(1- \pi_t^*
\bigr) \,dt, \qquad \pi_0^* = \pi_0
\end{equation}
[so $\pi_t^* = \pi_t^{h^*}$, where $(\pi_t^{h^*})$ is defined in
\eqref
{epiht1} with $h$ there replaced by $h^*$]. It is straightforward to
check that $(h^*, \tau^*)$ is indeed an optimal strategy (both in the
weak and strong formulations), and we do this in Section~\ref{sec5} in
the proof of Theorem~\ref{thmvf}.

On the other hand, in the case where $b < \lambda a \rho^2 / (2(a +
\lambda)^2)$, the optimal strategy should take the form mentioned in
\eqref{e1.4c},
%
%
\begin{equation}
\label{en1a0} h^*_t = 1_{\{\pi^*_t > A\}} \quad\mbox{and}\quad \tau^* = \inf
\bigl\{ t \geq0\dvtx\pi^*_t \geq B \bigr\},
\end{equation}
where the law of $(\pi^*_t)$ should be determined by the diffusion equation
%
%
\begin{equation}
\label{e6.69} d\pi^*_t = \lambda \bigl(1-\pi^*_t
\bigr) \,dt + \rho\pi^*_t \bigl(1-\pi^*_t \bigr)
1_{\{\pi
^*_t > A\}} \,d\bar W_t,
\end{equation}
or, looking back to \eqref{e6.69ab0} and \eqref{epiht1},
%
%
\begin{eqnarray}\label{e6.69ab}
d\pi^*_t &=& \lambda \bigl(1-\pi^*_t \bigr)
\,dt + \frac{r}{\sigma^2} \pi^*_t \bigl(1-\pi^*_t \bigr)
1_{\{\pi^*_t > A\}} (r 1_{\{\theta\leq t \}} \,dt + \sigma\,dW_t)
\nonumber
\\[-8pt]
\\[-8pt]
\nonumber
&&{} - \frac{r^2}{\sigma^2} \bigl(\pi^*_t \bigr)^2 \bigl(1-
\pi^*_t \bigr) 1_{\{\pi^*_t > A\}} \,dt.
\end{eqnarray}
Because of the irregularity of $p\mapsto1_{\{p > A\}}$, equations such
as \eqref{e6.69} and \eqref{e6.69ab} do not have a strong solution in
general (see, e.g., \cite{Chi,warren,KarShir}), but according
to the theory developed in \cite{GS}, Chapter~5, Section~24, they do
have a weak solution [such that the process $(\pi^*_t)$ spends an
amount of time at $A$ that has positive Lebesgue measure]. Therefore,
from the discussion in \eqref{e6.69ab0}--\eqref{tau0}, we expect
\eqref
{en1a0} to determine an optimal control system in the weak formulation
of our problem, but there will be no optimal strategy in the strong
formulation! This means that we will be able to use verification
Lemma~\ref{lem8} to prove, in Section~\ref{sec5}, that the function $g$
defined in Proposition~\ref{lem10} is equal to the value function
$\tilde g^{\weak}$, but a different approach via $\varepsilon$-optimal
strategies will be used to show that $g$ is equal to $\tilde g$.

\section{The value function}\label{sec5}

Formulas \eqref{defg} and \eqref{defg2} provide candidates, denoted by
$g$, for the value functions $\tilde g$ and $\tilde g^{\weak}$ defined,
respectively, in \eqref{vf1} and~\eqref{vf1a}. The objective of this
section is to prove that indeed, these two value functions are equal,
and equal to $g$.

%
\begin{thmm}\label{thmvf} \textup{(a)} \emph{Case where $0< b < \lambda a \rho
^2 / (2(a +
\lambda)^2)$.} Define $A$ by\break \eqref{e1.17}, let $f$ be as in \eqref
{e1.21a}, $K_1$ as in \eqref{e1.23}, $f_2$ as in \eqref{e1.22}, $B$ as
in \eqref{e1.24}, $K_2$ as in \eqref{e1.26}, $g_1$ as in \eqref{e1.28a}
and $g_2$ as in \eqref{e1.26a}. Then the function $g$ defined in
\eqref
{defg} is equal to the value function $\tilde g^{\weak}$ defined in
\eqref{vf1a}. Further, the control system associated to $h(t,p) = 1_{\{
p > A\}}$ and to $\tau^*$ in \eqref{en1a0} is optimal.

\textup{(b)} \emph{Case where $b \geq\lambda a \rho^2 / (2(a + \lambda)^2)$.}
Define $B$ by \eqref{e1.37a0} and $g_1$ by \eqref{e1.37a}. Then the
function $g$ defined in \eqref{defg2} is equal to the value function
$\tilde g^{\weak}$ defined in~\eqref{vf1a}.
\end{thmm}

%
\begin{thmm}\label{thmvf1} In both cases of Theorem~\ref{thmvf}, the
two value
functions $\tilde g$ (strong formulation) and $\tilde g^{\weak}$ (weak
formulation), defined, respectively, in \eqref{vf1} and~\eqref{vf1a}, are
equal (and equal to the function $g$ of Theorem~\ref{thmvf}).
\end{thmm}

%
\begin{remark}It is interesting to observe how the value function
$\tilde g$ and the thresholds $A$ and $B$ depend on the observation
cost $b$: we write $\tilde g(x,b)$, $A(b)$ and $B(b)$ to indicate this
dependence.

From \eqref{cost} and \eqref{vf1}, $b \mapsto\tilde g(x,b)$ is
nondecreasing. For $b=0$, $g(\cdot,0)$ is the value function obtained
in \cite{shiryaev}, Chapter~4.4, Theorem~9. As $b$ increases from $0$
to $b_c = \lambda a \rho^2/(2(a+b)^2)$, $B(b)$ decreases from $B(0)$ to
$B(b_c)$, and $A(b)$ increases from $0$ to $A(b_c) = \lambda
/(a+\lambda
) = B(b_c)$ [see \eqref{e1.17} for the first equality and the second
follows from the lines preceding \eqref{e1.24} since $f_2(A(b_c))
=-1$]. For $b \geq b_c$, $\tilde g(\cdot, b) = \tilde g(\cdot, b_c)$
since there is no dependence on $b$.
\end{remark}

Theorem~\ref{thmvf} will be proved in two steps. We begin by showing
that $g \leq\tilde g^{\weak}$.

%
\begin{lemma}\label{lemineq1} In both cases \textup{(a)} and \textup{(b)} of
Theorem~\ref{thmvf}, the
inequality $g \leq\tilde g^{\weak}$ holds.
\end{lemma}

\begin{pf} We are going to use part (1) of Lemma~\ref{lem8}. Suppose
first that we are in case (a) of Theorem~\ref{thmvf}. By construction,
and in particular by \eqref{e1.10}, \eqref{e1.13} and~\eqref{e1.13a},
$g$ is $C^2$ on $[0,B[$, and $C^1$
on $[0,1]$ by \eqref{e1.12} and \eqref{e1.7}, so
%
\begin{equation}
\label{e5.65} g'(B-) = g'(B+) = -1.
\end{equation}
%
By \eqref{e1.29a}, $0\leq g(x) \leq1-x$. Let $((h_t), \tau, X)$ be a
control system, and set
%
%
\begin{equation}
\label{e5.65aa} Y_t = g \bigl(\pi^h_t \bigr)
+ a \int_0^t \pi^h_s
\,ds + b \int_0^t h_s \,ds.
\end{equation}
We now apply It\^o's formula, in the form given in \cite{PS}, Section~3.5:
%
%
\begin{eqnarray}
\label{e5.66aa}
Y_t &=& Y_0 + \int
_0^t g' \bigl(
\pi^h_s \bigr) \,d\pi^h_s + \int
_0^t \bigl(a \pi^h_s +
b h_s \bigr) \,ds + \frac{1}{2}\int_0^t
g'' \bigl(\pi^h_s \bigr) \,d
\bigl\langle\pi^h \bigr\rangle_s
\nonumber
\\[-8pt]
\\[-8pt]
\nonumber
&&{} + \frac{1}{2} \bigl(g'(B+) - g'(B-) \bigr)
L^B_t,
\end{eqnarray}
where $ L^B_t$ is the local time of $(\pi^h_s)$ at $B$. By \eqref
{e5.65}, the factor $g'(B+) - g'(B-)$ vanishes, so as in \eqref{drift},
we find that
%
%
\begin{equation}
\label{e5.66} Y_t = Y_0 + \int_0^t
g' \bigl(\pi^h_s \bigr) \frac{r}{\sigma}
\pi^h_s \bigl(1- \pi^h_s \bigr)
\sqrt{h_s} \,d\bar W_s + \int_0^t
\Phi \bigl(\pi^h_s, h_s \bigr) \,ds,
\end{equation}
where
\[
\Phi(x,\eta) = L_1 g(x) + \eta \bigl[\tfrac{1}{2}g''(x)
\bigl(\rho x (1-x) \bigr)^2 + b \bigr],\qquad \eta\in[0,1],
\]
and $L_1$ is defined in \eqref{e1.28b}. We note that by construction
and by Proposition~\ref{lem10},
\begin{eqnarray*}
\Phi(x,0) &=& L_1 g(x) \geq0 \qquad\mbox{for all } x \in[0,1],
\\
\Phi(x,1) &=& L_2 g(x) \geq0\qquad \mbox{for all } x \in[0,1]\setminus\{
B\},
\end{eqnarray*}
where $L_2$ is defined in \eqref{e1.28c}, and since $\eta\mapsto\Phi
(x,\eta)$ is an affine function, we conclude that $\Phi(x,\eta) \geq
0$, for all $\eta\in[0,1]$. Since $g'$ is bounded on $[0,1]$, the
stochastic integral in \eqref{e5.66} is an $\F^X$-martingale [recall
\eqref{emt}], and therefore $(Y_t)$ is an $\F^X$-submartingale. The
conclusion now follows from part (1) of Lemma~\ref{lem8}.

Now suppose that we are in case (b) of Theorem~\ref{thmvf}. By
construction, $g$ is $C^2$ on $[0,B[$, and $C^1$ on $[0,1]$ by \eqref
{e1.35b} and \eqref{e1.35c}, so
\[
g'(B-) = g'(B+) = -1.
\]
By \eqref{e1.38a}, $0 \leq g(x) \leq1-x$, for all $x \in[0,1]$. Let
$((h_t),\tau, X)$ be a control system, and define $Y_t$ as in \eqref
{e5.65aa}. Applying It\^o's formula, we obtain \eqref{e5.66aa}, and
this leads again to \eqref{e5.66}. Using this time Proposition~\ref
{lem11}, we see that $\Phi(x,\eta) \geq0$, for all $\eta\in[0,1]$.
Therefore, we conclude, as before, that $(Y_t)$ is an $\F
^X$-submartingale, and the conclusion follows from part (1) of
Lemma~\ref{lem8}.
\end{pf}

We now prove Theorem~\ref{thmvf}.

\begin{pf*}{Proof of Theorem~\ref{thmvf}} We begin with case (b).
As mentioned in \eqref{en1a} and \eqref{en1}, the candidate optimal
control system is $(h^*,\tau^*, X^*)$, where $h_t^* \equiv0$, $X^*
\equiv0$ and
\[
\tau^* = \inf \bigl\{ t \geq0\dvtx\pi_t^* \geq B \bigr\},
\]
where $(\pi_t^*)$ is defined in \eqref{en1}.
Clearly, $(h^*, \tau^*, X^*)$ is a control system, and so it suffices
to check properties (b), (c) and (d) of Lemma~\ref{lem8}. By \eqref{e1.38aa},
\[
dY^*_t = L_1 g \bigl(\pi_t^* \bigr) \,dt =
0 \qquad\mbox{for } t < \tau^*.
\]
Therefore, $(Y^*_{t\wedge\tau^*})$ is a (constant and deterministic)
martingale, proving (b).

Further, since $(\pi_t^*)$ is deterministic, we solve \eqref{en1} to
find that
%
%
\begin{equation}
\label{eln} \tau^* = \cases{ %
\displaystyle\frac{1}{\lambda} \ln \biggl(\frac{1-\pi_0}{1-B} \biggr), &\quad
$\mbox {if } \pi
_0 < B,$
\vspace*{2pt}\cr
0, &\quad  $\mbox{if } \pi_0 \geq B,$}
\end{equation}
so (c) holds. Finally, if $\pi_0 < B$, then
\[
g \bigl(\pi^*_{\tau^*} \bigr) = g(B) = 1-B = 1-\pi^*_{\tau^*}
\]
by \eqref{e1.35b}, and if $\pi_0 \geq B$, then
\[
g \bigl(\pi^*_{\tau^*} \bigr) = g \bigl(\pi^*_{0} \bigr) = 1-
\pi^*_{0} = 1- \pi^*_{\tau^*}
\]
by \eqref{defg2}. This proves case (b) of Theorem~\ref{thmvf}.


We now consider case (a). We have seen in Lemma~\ref{lemineq1} that $g
\leq\tilde g^{\weak}$. In order to establish the converse inequality,
consider $(h^*_t)$ and $\tau^*$ defined in \eqref{en1a0} and the
associated control system $((h^*_t),\tau^*,X^*)$, constructed as in
\eqref{e6.69ab0}--\eqref{tau0}, using the function $h(t,p) = 1_{\{p >
A\}}$ and $\pi^*_t$ defined as a weak solution of \eqref{e6.69ab}. Then
for $t \leq\tau^*$,
\begin{eqnarray*}
\Phi \bigl(\pi^*_t,h^*_t \bigr) &=& \cases{
L_1 g \bigl(
\pi^*_t \bigr), & \quad$\mbox{if } \pi^*_t < A,$
\vspace*{2pt}\cr
L_2 g \bigl(\pi^*_t \bigr), &\quad $\mbox{if }
\pi^*_t \in[A,B[,$}
\\
& =& 0
\end{eqnarray*}
by \eqref{e1.29aa} and \eqref{e1.29aaa}. Therefore, $(Y^*_{t \wedge
\tau
^*})$ is an $\F^X$-martingale. According to Lemma~\ref{lemtaustar}
below, $E_\pi(\tau^*) < \infty$, and $g(\pi^{h^*}_{\tau^*}) = g(B) =
1-B$ by \eqref{en1a0} and \eqref{defg}. This proves properties (b), (c)
and (d) of Lemma~\ref{lem8} and concludes the proof that $ g = \tilde
g^{\weak}$ and $((h^*_t),\tau^*,X^*)$ is an optimal control system,
since we already verified~(a) of Lemma~\ref{lem8} during the proof of
Lemma~\ref{lemineq1}.
\end{pf*}

%
\begin{lemma}\label{lemtaustar} Suppose that we are in case \textup{(a)} of
Theorem~\ref{thmvf}.
Let $\tau^*$ be defined as in \eqref{en1a0}. Then for all $\pi\in
[0,1]$, $E_\pi(\tau^*) < \infty$.
\end{lemma}

\begin{pf} If $\pi\in[B,1]$, then $E_\pi(\tau^*) = 0$, and if $\pi
\in[0,A[$, then $(\pi^*_t)$ reaches $A$ at the deterministic time
$\lambda^{-1} \ln((1-\pi)/(1-A))$ [see \eqref{eln}], so the problem
reduces to considering $\pi\in[A,B[$.

Recall from \eqref{e3.14} and \eqref{e6.69} that $(\pi^*_t)$ solves,
in the terminology of \cite{GS}, Chapter~5, Section~24, an s.d.e. with
\emph{delayed} reflection at the boundary point $A$, and this process
is associated to a diffusion $(\tilde\xi_t)$ with \emph{instantaneous}
reflection at the boundary
%
%
\begin{equation}
\label{exi5.6} d\tilde\xi_t = \lambda(1-\tilde\xi_t)
\,dt + \rho\tilde\xi_t (1-\tilde\xi_t) \,d\tilde
W_t + d \zeta_t,
\end{equation}
where $(\zeta_t)$ is a nondecreasing process that increases at those
points where $\tilde\xi_t = A$,
\[
\tilde W_t = \int_0^{\phi_t}
1_{\{\tilde\xi_{\tau_s} > A \}} \,d \bar W_s,
\]
and $\phi_t$ is defined by the relation
\[
t = \int_0^{\phi_t} 1_{\{\tilde\xi_{\tau_s} > A \}} \,ds.
\]
As explained in \cite{GS}, $(\pi^*_t)$ has the same law as $(\tilde
\xi
_{\tau_t})$, where $\tau_t$ is defined by the relation
\[
t = \tau_t + \frac{1}{\lambda(1-A)} \zeta_{\tau_t}.
\]
Therefore, $\tau^*$ has the same law as $T = \inf\{t \in\IR_+\dvtx
\tilde
\xi_{\tau_t} = B\}$. Letting $\sigma= \inf\{s \in\IR_+\dvtx
\tilde\xi_s
= B\}$, we see that $\sigma= \tau_T$. Further, according to Lemma~5 in
\cite{GS}, Chapter~5, Section~23, $E_\pi(\sigma) = V_0(y) - V_0(\pi)$,
where $V_0'(A) = 0$ and for $y \in\,]A,1[$,
\[
\lambda(1-y) V'_0(y) + \tfrac{1}{2}
\rho^2 y^2 (1-y)^2 V_0''(y)
= 1.
\]
An explicit expression for $V_0$ can be obtained by using \eqref{e1.22}
with $K_1 = 0$, $a=0$ and $b=-1$, and then integrating from $A$ to $y$.
In particular, $E_\pi(\sigma) < +\infty$.

Notice that
\[
T = \tau_T + \frac{1}{\lambda(1-A)} \zeta_{\tau_T} = \sigma+
\frac
{1}{\lambda(1-A)} \zeta_\sigma.
\]
Therefore, it suffices to show that $E_\pi(\zeta_\sigma) < +\infty
$. By
\eqref{exi5.6},
%
%
\begin{equation}
\label{xitmin} \tilde\xi_{t\wedge\sigma} = \tilde\xi_0 + \int
_0^{t\wedge\sigma} \lambda(1-\tilde\xi_s) \,ds
+ \int_0^{t\wedge\sigma} \rho\tilde\xi_s (1-
\tilde\xi_s) \,d\tilde W_s + \zeta_{t\wedge\sigma}.
\end{equation}
The stochastic integral is an $L^2$-bounded martingale, since
\[
E_\pi \biggl( \int_0^{t\wedge\sigma}
\rho^2 \tilde\xi_s^2 (1-\tilde\xi
_s)^2 \,ds \biggr) \leq\rho^2
E_\pi(\sigma) < +\infty.
\]
Therefore, the optional sampling theorem can be applied and, since the
$ds$-integral in \eqref{xitmin} is nonnegative, we find that
\[
B = E_\pi(\tilde\xi_\sigma) \geq\pi+ E_\pi(
\zeta_\sigma),
\]
so $E_\pi(\zeta_\sigma) < +\infty$, as was to be proved.
\end{pf}

For the remainder of this section, we put ourselves in case (a) of
Theorem~\ref{thmvf}. Since we have observed just after \eqref{vf1a}
that $\tilde g \geq\tilde g^{\weak}$, and $\tilde g^{\weak} = g$ by
Theorem~\ref{thmvf}(a), in order to prove Theorem~\ref{thmvf1}, it
suffices to establish the inequality $g \geq\tilde g$. For $\ep>0$,
we are going to define an admissible control $h^\ep$, and a strategy
$(h^\ep, \tau^\ep)$, with associated cost $\tilde g_\ep= E(C(h^\ep,
\tau^\ep))$, and we shall show that $\tilde g_\ep\to g$ as $\ep
\downarrow0$. From the definition of $\tilde g$ in \eqref{vf1}, this
will establish that $g \geq\tilde g$, and this will prove Theorem~\ref
{thmvf1}.

\subsection*{An almost optimal strategy}

Define the function
\[
h^{(\ep)}(x) = \frac{x-A}{\ep} 1_{]A,A+\ep[}(x) +
1_{[A+\ep, \infty[}(x).
\]
%

Consider the s.d.e. 
%
\begin{eqnarray}
\label{e5.67}
dp^\ep_t &=& \lambda \bigl(1-
p^\ep_t \bigr) \,dt\nonumber\\
&&{}+ \frac{r}{\sigma^2}
p^\ep_t \bigl(1-p^\ep_t \bigr)
\Bigl(r h^{(\ep)} \bigl(p^\ep_t \bigr)
1_{\{\theta< t\}} \,dt + \sigma\sqrt {h^{(\ep)}
\bigl(p^\ep_t \bigr)} \,dW_t \Bigr)
\\
&&{} - \frac{r^2}{\sigma^2} \bigl(p^\ep_t \bigr)^2
\bigl(1- p^\ep_t \bigr) h^{(\ep)}
\bigl(p^\ep_t \bigr) \,dt,\nonumber
\end{eqnarray}
with $p^\ep_0 = \pi_0$. According to \cite{IW}, Theorem~3.2 page 168,
this s.d.e. has a unique strong solution $(p^\ep_t, t \geq0)$, since
$\sqrt{h^{(\ep)}}$ is H\"older-continuous with exponent $1/2$.

Set
%
%
\begin{equation}
\label{e5.70} \tau^\ep= \cases{ %
\inf \bigl\{t \geq0\dvtx p^\ep_t \geq B
\bigr\}, & \quad$\mbox{if } \{ \cdots\} \neq\varnothing,$
\vspace*{2pt}\cr
+\infty,& \quad$\mbox{otherwise.}$}
\end{equation}
Using \eqref{e6.69ab0}--\eqref{tau0}, we associate to $(h^{(\ep
)},\tau
^\ep)$ a strategy $((h^\ep_t),\tau^\ep)$.

We are now going to determine the cost of the strategy $(h^\ep, \tau
^\ep)$, and we will see in Proposition~\ref{lem15} below that for
$\ep$
small, this strategy is nearly optimal. Let
%
%
\begin{equation}
\label{e5.71} \tilde g_\ep(\pi_0) = E \bigl(C
\bigl(h^\ep, \tau^\ep \bigr) \bigr).
\end{equation}
In order to determine the function $\tilde g_\ep$, we will use the
following lemma.

%
\begin{lemma}\label{lem14} Suppose that we are in case \textup{(a)} of
Theorem~\ref{thmvf}
and that we can find a continuous function $g_\ep$ on $[0,1]$ that is
$C^2$ on $[0,1]\setminus\{A,B\}$, $C^1$ on $[0,1]\setminus\{B\}$ and
such that
%
%
\begin{equation}
\label{e5.72} Lg_\ep(x) = - \bigl(a x + b h^{(\ep)}(x)
\bigr),
\end{equation}
where $Lg_\ep(x)$ is defined by
%
%
\begin{equation}
\label{e5.73} Lg_\ep(x) = \lambda(1-x) g'_\ep(x)
+ \tfrac{1}{2}\rho^2 x^2 (1-x)^2
h^{(\ep
)}(x) g''_\ep(x)
\end{equation}
and
%
%
\begin{equation}
\label{e5.74} g_\ep(x) = 1-x\qquad \mbox{for } x \in[B,1].
\end{equation}
If, in addition,
%
%
\begin{equation}
\label{e5.75} E_x(\tau_\ep) < +\infty\qquad \mbox{for all }
x \in[0,1],
\end{equation}
then $g_\ep= \tilde g_\ep$.
\end{lemma}

\begin{pf} Suppose $\pi_0 \in[B,1]$. Then $g_\ep(\pi_0) = 1 - \pi_0$,
and since $\tau^\ep=0$ a.s., \eqref{cost} gives
\[
\tilde g_\ep(\pi_0) = E \bigl(C \bigl(h^\ep,
\tau^\ep \bigr) \bigr) = P\{\theta>0\} = 1-\pi_0.
\]
Therefore, by \eqref{e5.74}, $g_\ep(\pi_0) = \tilde g_\ep(\pi_0)$ in
this case.

Now suppose that $\pi_0 \in[0,B[$. According to Lemma~\ref{lem0} and
\eqref{ept},
\[
\tilde g_\ep \bigl(p^\ep_0 \bigr) = E \bigl(C
\bigl(h^\ep, \tau^\ep \bigr) \bigr) = E \biggl(1 -
p^{\ep
}_{\tau^\ep} + a \int_0^{\tau^\ep}
p^{\ep}_s \,ds + b \int_0^{\tau^\ep}
h^{(\ep)} \bigl(p^{\ep}_s \bigr) \,ds \biggr).
\]
Since $\tau^\ep< +\infty$ a.s. by \eqref{e5.75}, $p^{\ep}_{\tau
^\ep} =
B$, and $1-B = g_\ep(B)$ by \eqref{e5.74}, so
\[
E \bigl(C \bigl(h^\ep, \tau^\ep \bigr) \bigr) = E \biggl(
g_\ep \bigl(p^\ep_{\tau^\ep} \bigr) + a \int
_0^{\tau
^\ep} p^{\ep}_s \,ds + b
\int_0^{\tau^\ep} h^{(\ep)}
\bigl(p^{\ep}_s \bigr) \,ds \biggr).
\]

As in Lemmas \ref{lem6} and \ref{lem7}, we see from \eqref{e5.67} and
\eqref{e6.69ab2} that
%
%
\begin{equation}
\label{e5.76} d p^\ep_t = \lambda \bigl(1-p^\ep_t
\bigr) \,dt + \rho p^\ep_t \bigl(1-p^\ep_t
\bigr) \sqrt{h^{(\ep)} \bigl(p^\ep_t
\bigr)} \,d \bar W^\ep_t,
\end{equation}
where $(\bar W^\ep_t)$ is an Brownian motion. Let
\[
M_t = g_\ep \bigl(p^\ep_t \bigr) +
a \int_0^t p^\ep_s
\,ds + b \int_0^t h^{(\ep
)}
\bigl(p^\ep_s \bigr) \,ds.
\]
We apply It\^o's formula in the form given in \cite{PS}, Section~3.5.3,
using the fact that $g_\ep'(A-) = g_\ep'(A+)$, to see that
\begin{eqnarray*}
dM_t &=& g'_\ep \bigl(p^\ep_t
\bigr) \Bigl[\lambda \bigl(1 - p^\ep_t \bigr) \,dt + \rho
p^\ep_t \bigl(1 - p^\ep_t \bigr)
\sqrt{h^{(\ep)} \bigl(p^\ep_t \bigr)} \,d
\bar W^\ep_t \Bigr]
\\
&&{} + \bigl[\tfrac{1}{2} g''
\bigl(p^\ep_t \bigr) \rho^2
\bigl(p^\ep_t \bigr)^2 \bigl(1-p^\ep_t
\bigr)^2 h^{(\ep)} \bigl(p^\ep_t \bigr)
+ a p^\ep_t + b h^\ep \bigl(p^\ep_t
\bigr) \bigr] \,dt
\\
&=& g'_\ep \bigl(p^\ep_t \bigr)
\rho p^\ep_t \bigl(1 - p^\ep_t
\bigr) \sqrt{h^{(\ep)} \bigl(p^\ep_t
\bigr)} \,d\bar W^\ep_t + \bigl[Lg_\ep
\bigl(p^\ep_t \bigr) + a p^\ep_t + b
h^\ep \bigl(p^\ep_t \bigr) \bigr] \,dt.
\end{eqnarray*}
By \eqref{e5.72}, the drift in brackets vanishes, and therefore $(M_{t
\wedge\tau^\ep}, t \geq0)$ is an $\F^{X^\ep}$-martingale. Since, by
\eqref{ept}, $0\leq p^\ep_s \leq1$, and $0 \leq h^{(\ep)} \leq1$ and
$g_\ep$ is bounded, we see that $\vert M_t\vert\leq A + (a+b)t$, so
$\vert M_{t\wedge\tau^\ep} \vert\leq A + (a+b)\tau^\ep$. Since
$E_{\pi
_0}(\tau^\ep) < +\infty$ by \eqref{e5.75}, $(M_{t\wedge\tau^\ep
})$ is
uniformly integrable, and so
\[
\tilde g_\ep(\pi_0) = E \bigl(C \bigl(h^\ep,
\tau^\ep \bigr) \bigr) = E_{\pi_0}(M_{\tau^\ep}) =
E_{\pi_0}(M_0) = g_\ep \bigl(p^\ep_0
\bigr) = g_\ep(\pi_0).
\]
Therefore, $g_\ep(\pi_0) = \tilde g_\ep(\pi_0)$ as claimed. This
completes the proof of Lemma~\ref{lem14}.
\end{pf}

\subsection*{Constructing $g_\ep$}

It remains to construct the function $g_\ep$ satisfying the
assumptions of Lemma~\ref{lem14}. Notice that on $]0,A[$, writing
$\bar
g_1$ instead of $g_\ep$, equation \eqref{e5.72} becomes
%
%
\begin{equation}
\label{e5.77} \lambda(1-x) \bar g'_1(x) + ax = 0,
\end{equation}
and as in \eqref{e1.37}, the solution of this differential equation is
%
%
\begin{equation}
\label{e5.78} \bar g_1(x) = \frac{a}{\lambda} \bigl(x + \ln(1-x)
\bigr) + K^\ep_1,
\end{equation}
where $K^\ep_1$ is a constant to be determined.

On $]A+\ep, B[$, writing $\bar g_3$ instead of $g_\ep$, equation
\eqref
{e5.72} becomes
%
%
\begin{equation}
\label{e5.79} \lambda(1-x) \bar g'_3(x) +
\tfrac{1}{2} \rho^2 x^2 (1-x)^2 \bar
g''_3(x) + a x + b = 0,
\end{equation}
which is the same equation as in \eqref{e1.9}, and as in \eqref{e1.25},
its solution is
%
%
\begin{equation}
\label{e5.80} \bar g_3(x) = \int_{A + \ep}^x
\bar h_3(y) \,dy + K_3^\ep,
\end{equation}
where
%
%
\begin{equation}
\label{e5.81} \bar h_3(x) = K^\ep_2 f(x) +
f(x) \int_{A+\ep}^x \frac{-2}{\rho^2}
\frac{ay + b}{y^2(1-y)^2} \frac{1}{f(y)} \,dy
\end{equation}
and
%
%
\begin{equation}
\label{e5.82} f(x) = \biggl(\frac{1-x}{x} \biggr)^\alpha
e^{\alpha/x}\qquad \mbox{where } \alpha= \frac{2 \lambda}{\rho^2},
\end{equation}
and $K^\ep_2$, $K^\ep_3$ are constants to be determined.

Finally, on $]A,A+\ep[$, writing $\bar g_2$ instead of $g_\ep$,
equation \eqref{e5.72} becomes
%
%
\begin{equation}\quad
\label{e5.83} \lambda(1-x) \bar g'_2(x) +
\frac{1}{2} \rho^2 x^2 (1-x)^2
\frac{x-A}{\ep} \bar g''_2(x) + ax + b
\frac{x-A}{\ep} = 0.
\end{equation}
Let $\bar h_2(x) = \bar g'_2(x)$, so the associated homogeneous
equation is
%
%
\begin{equation}
\lambda(1-x) \bar f_2(x) + \frac{1}{2}\rho^2
x^2 (1-x)^2 \frac{x-A}{\ep} \bar f'_2(x)
= 0,
\end{equation}
whose solution is
%
%
\begin{equation}
\label{e5.85} \bar f^\ep_2(x) = \psi_\ep(x)
(x-A)^{- \beta_\ep},
\end{equation}
where
\[
\beta_\ep= \frac{1}{A^2 (1-A)} \frac{2 \lambda\ep}{\rho^2}
\]
and
\[
\psi_\ep(x) = x^{2 \lambda\ep(1+A)/(\rho A)^2} (1-x)^{2 \lambda\ep
/(\rho^2 (1-A))} \exp \biggl(-
\frac{2 \lambda\ep}{A \rho^2} \frac{1}{x} \biggr).
\]
Therefore,
%
%
\begin{eqnarray}\qquad
\label{e5.85a} \bar h_2(x)& = &K \bar f^\ep_2(x)
\nonumber
\\[-8pt]
\\[-8pt]
\nonumber
&&{}+ \bar f^\ep_2(x) \int_{A}^x
\frac
{-2\ep}{\rho^2} \biggl(ay + b \frac{y-A}{\ep} \biggr) \frac{1}{y^2(1-y)^2
(y-A)}
\frac{1}{\bar f^\ep_2(y)} \,dy,
\end{eqnarray}
and if we want $\bar h_2$ to be bounded as $x \downarrow A$, then we
must set $K=0$ (notice that there is no integrability problem at
$y=A$). We conclude that
%
%
\begin{equation}\qquad
\label{e5.86} \bar h_2(x) = \bar f^\ep_2(x)
\int_{A}^x \frac{-2\ep}{\rho^2} \biggl(ay + b
\frac{y-A}{\ep} \biggr) \frac{1}{y^2(1-y)^2 (y-A)} \frac{1}{\bar
f^\ep_2(y)} \,dy
\end{equation}
and
%
%
\begin{equation}
\label{e5.87} \bar g_2(x) = \int_A^x
\bar h_2(y) \,dy + K^\ep_4,
\end{equation}
where $K^\ep_4$ is a constant to be determined.

In order to determine the four constants $K^\ep_1,\dots,K^\ep_4$, we
shall impose the four equations
%
%
\begin{eqnarray}
\label{e5.88} \bar g_3 (B) &=& 1-B,
\\
\label{e5.89} \bar g_3(A+\ep) &=& \bar g_2(A+\ep),
\\
\label{e5.90} \bar g'_3(A+\ep) &=& \bar
g'_2(A+\ep),
\\
\label{e5.91} \bar g_1(A) &=& \bar g_2(A).
\end{eqnarray}
We note that \eqref{e5.89} and \eqref{e5.90}, together with \eqref
{e5.79} and \eqref{e5.83}, imply that $\bar g''_3(A+\ep) = \bar
g''_2(A+\ep)$, and \eqref{e5.91}, together with \eqref{e5.77} and
\eqref
{e5.83}, implies that $\bar g'_1(A) = \bar g'_2(A)$.

From \eqref{e5.88} and \eqref{e5.80}, we see that
%
%
\begin{equation}
\label{e5.92} K^\ep_3 = 1 - B - \int
_{A+\ep}^B \bar h_3(y) \,dy,
\end{equation}
while \eqref{e5.89}, \eqref{e5.80} and \eqref{e5.87} imply that
%
%
\begin{equation}
\label{e5.93} K^\ep_3 = \int_A^{A+\ep}
\bar h_2(y) \,dy + K^\ep_4.
\end{equation}
Equality \eqref{e5.90}, \eqref{e5.81} and \eqref{e5.86} give the relation
%
%
\begin{eqnarray}
\label{e5.94} K^\ep_2& =& \frac{\bar f^\ep_2(A+\ep)}{f(A+\ep)}
\nonumber
\\[-8pt]
\\[-8pt]
\nonumber
&&{}\times\int
_{A}^{A+\ep} \frac
{-2\ep}{\rho^2} \biggl(ay + b
\frac{y-A}{\ep} \biggr) \frac{1}{y^2(1-y)^2
(y-A)} \frac{1}{\bar f^\ep_2(y)} \,dy,
\end{eqnarray}
while \eqref{e5.91}, \eqref{e5.78} and \eqref{e5.87} give
%
%
\begin{equation}
\label{e5.95} \frac{a}{\lambda} \bigl(A + \ln(1-A) \bigr) +
K^\ep_1 = K^\ep_4.
\end{equation}
Therefore, \eqref{e5.94} determines $K^\ep_2$, \eqref{e5.92} determines
$K^\ep_3$, then \eqref{e5.93} determines $K^\ep_4$ and \eqref{e5.95}
determines $K^\ep_1$.

%
\begin{prop}\label{lem15} For $\ep>0$, let $K^\ep_1,\dots,K^\ep_4$
be determined by
\eqref{e5.92}--\eqref{e5.95}, define $\bar g_1(x)$ as in \eqref{e5.78},
$\bar g_2(x)$ as in \eqref{e5.87}, and $\bar g_3(x)$ as in \eqref
{e5.80}. Set
\[
g_\ep(x) = \cases{ %
\bar
g_1(x), & \quad$\mbox{if } 0 \leq x \leq A,$
\vspace*{2pt}\cr
\bar g_2(x), &\quad $\mbox{if } A < x < A+\ep,$
\vspace*{2pt}\cr
\bar g_3(x), &\quad $\mbox{if } A+\ep\leq x < B,$
\vspace*{2pt}\cr
1-x, & \quad$\mbox{if } B \leq x \leq1.$}
\]
Then $g_\ep$ satisfies the assumptions of Lemma~\ref{lem14}. Further,
let $g$ be as in case \textup{(a)} of Theorem~\ref{thmvf}. Then
\[
\lim_{\ep\downarrow0} g_\ep(x) = g(x)\qquad\mbox{for all } x
\in[0,1].
\]
\end{prop}

\begin{pf} By the comments that follow \eqref{e5.91}, $g_\ep$ is $C^2$
on $[0,1]\setminus\{A,B\}$, $C^1$ on $[0,1]\setminus\{B\}$ and
continuous on $[0,1]$. For $x \in[B,1]$, $g_\ep(x) = 1-x = g(x)$, so
we consider the case where $x \in[0,B[$.

\emph{Case} 1: $x \in\,]A,B[$. We first check that $K^\ep_2
\to K_1$, where $K_1$ is defined in \eqref{e1.23}. We note that
\begin{eqnarray*}
K^\ep_2 &=& \frac{\psi_\ep(A+\ep)}{f(A+\ep)} \ep^{-\beta_\ep} \\
&&{}\times\int
_{A}^{A+\ep} \frac{-2\ep}{\rho^2} \biggl(ay + b
\frac{y-A}{\ep} \biggr) \frac{1}{y^2(1-y)^2} \frac{1}{\psi_\ep(y)}
(y-A)^{\beta_\ep- 1} \,dy.
\end{eqnarray*}
Notice that $\psi_\ep(A+\ep) \to1$ and $f(A+\ep) \to f(A)$ as $\ep
\downarrow0$. Set
\[
\lambda_0 = \frac{1}{A^2(1-A)} \frac{2\lambda}{\rho^2}\qquad\mbox{so that }
\beta_\ep= \lambda_0 \ep.
\]
Then
\begin{eqnarray*}
K_2^\ep&\sim&\frac{1}{f(A)} \ep^{1-\lambda_0 \ep} \int
_{A}^{A+\ep} \frac{-2}{\rho^2} \biggl(ay + b
\frac{y-A}{\ep} \biggr) \frac
{1}{y^2(1-y)^2 \psi_\ep(y) } (y-A)^{\lambda_0 \ep- 1} \,dy
\\
&\sim&\frac{1}{f(A)} \frac{-2}{\rho^2} \frac{1}{A^2(1-A)^2 \psi
_\ep
(A)} \\
&&{}\times\ep^{1-\lambda_0 \ep}
\int_{A}^{A+\ep} \biggl[a A (y-A)^{\lambda_0
\ep-1} +
\frac{b}{\ep} (y-A)^{\lambda_0 \ep} \biggr] \,dy,
\end{eqnarray*}
and the integral is equal to
\[
aA \frac{\ep^{\lambda_0 \ep}}{\lambda_0 \ep} + \frac{b}{\ep} \frac{\ep
^{\lambda_0 \ep+1}}{\lambda_0 \ep+1},
\]
and therefore,
\[
\lim_{\ep\downarrow0} K^\ep_2 =
\frac{1}{f(A)} \frac{-2}{\rho^2} \frac
{1}{A^2(1-A)^2} \frac{aA}{\lambda_0} = -
\frac{a}{\lambda} \frac
{A}{1-A} \frac{1}{f(A)} = K_1,
\]
as claimed.

This implies that for $y > A + \ep$, $\bar h_3(y) \to f_2(y)$, where
$\bar h_3$ and $f_2$ are, respectively, defined in \eqref{e5.81} and
\eqref{e1.22}. By dominated convergence, we deduce that $K^\ep_3 \to
K_2$, and for $x \in\,]A,B[$ and for $\ep\downarrow0$ with $0<\ep< x-A$,
\[
g_\ep(x) = \bar g_3(x) \to g_2(x) = g(x),
\]
where $g_2$ is defined in \eqref{e1.26a}.


\emph{Case} 2: $x \in[0,A]$. 
From \eqref{e5.93}, we see that $K^\ep_3 - K^\ep_4 \to0$, therefore
$K^\ep_4 \to K_2$ by the above, and from \eqref{e5.95}, we see that
\[
K^\ep_1 \to K_2 - \frac{a}{\lambda} \bigl(A
+ \ln(1-A) \bigr).
\]
We conclude from \eqref{e5.78} and \eqref{e1.28a} that for $x \in
[0,A]$, as $\ep\downarrow0$,
\[
g_\ep(x) = \bar g_1(x) \to g_1(x) = g(x).
\]
%
This completes the proof of Proposition~\ref{lem15}.
\end{pf}

The next lemma checks condition \eqref{e5.75}.

%
\begin{lemma}\label{lem16} Fix $\ep>0$, and let $\tau^\ep$ be
defined in \eqref
{e5.70}. Then for all $x \in[0,1]$, $E_x(\tau^\ep) < \infty$.
\end{lemma}

\begin{pf} We first seek a bounded function $\gamma_\ep$ defined on
$[0,B]$ such that
%
%
\begin{equation}
\label{e5.96} L \gamma_\ep= -1,
\end{equation}
where $L$ is the operator defined in \eqref{e5.73}.

For $0<x<A$, \eqref{e5.96} becomes
%
%
\begin{equation}
\label{e5.97} \lambda(1-x) \gamma'_\ep(x) = -1,
\end{equation}
so
%
%
\begin{equation}
\label{e5.98} \gamma_\ep(x) = \frac{1}{\lambda} \ln(1-x) +
D_1,\qquad 0\leq x \leq A.
\end{equation}

For $A < x < A+\ep$, \eqref{e5.96} becomes
%
%
\begin{equation}
\label{e5.99} \lambda(1-x) \gamma'_\ep(x) +
\frac{1}{2} \rho^2 x^2(1-x)^2
\frac{x-A}{\ep} \gamma''_\ep(x) = -1,
\end{equation}
and as in \eqref{e5.83} and \eqref{e5.85a}, the solution to this
equation is
%
%
\begin{equation}
\label{e5.100} \gamma_\ep(x) = \int_A^x
h_4(y) \,dy + D_3, \qquad A < x < A+\ep,
\end{equation}
where
%
%
\begin{equation}
\label{e5.101} h_4(x) = D_2 \bar f^\ep_2(x)
+ \bar f^\ep_2(x) \int_{A}^x
\frac{-2\ep
}{\rho^2} \frac{1}{y^2(1-y)^2 (y-A)} \frac{1}{\bar f^\ep_2(y)} \,dy,
\end{equation}
and $\bar f^\ep_2$ is defined in \eqref{e5.85}. Since we want $h_4$ and
$\gamma_\ep$ to be bounded (as $x \downarrow A$), we set $D_2 = 0$.

For $A+\ep< x < B$, \eqref{e5.96} becomes
%
%
\begin{equation}
\label{e5.102} \lambda(1-x) \gamma'_\ep(x) +
\tfrac{1}{2} \rho^2 x^2(1-x)^2
\gamma''_\ep(x) = -1,
\end{equation}
and as in \eqref{e5.80}, the solution of this equation is
%
%
\begin{equation}
\label{e5.103} \gamma_\ep(x) = \int_{A+\ep}^x
h_5(y) \,dy + D_4,
\end{equation}
where
%
%
\begin{equation}
\label{e5.104} h_5(x) = D_5 f(x) + f(x) \int
_{A+\ep}^x \frac{-2}{\rho^2 y^2(1-y)^2} \frac{1}{f(y)}
\,dy,
\end{equation}
and $f(x)$ is defined in \eqref{e5.82}.

We must determine the constants $D_1,\dots,D_5$. For this, we impose
the following conditions:

\begin{longlist}[(a)]
\item[(a)] $\gamma_\ep(B) = 0$,

\item[(b)] $\gamma_\ep((A+\ep)+) = \gamma_\ep((A+\ep)-)$,

\item[(c)] $\gamma'_\ep((A+\ep)+) = \gamma'_\ep((A+\ep)-)$,

\item[(d)] $\gamma_\ep(A+) = \gamma_\ep(A-)$.
\end{longlist}

We note that (b) and (c), together with \eqref{e5.99} and
\eqref{e5.102}, imply that
%
%
\begin{equation}
\gamma''_\ep \bigl((A+\ep)+ \bigr) =
\gamma''_\ep \bigl((A+\ep)- \bigr),
\end{equation}
so $\gamma_\ep$ will be $C^2$ at $A+\ep$. Also, (d) together with
\eqref
{e5.97} and \eqref{e5.99} implies that
%
%
\begin{equation}
\gamma'_\ep(A+) = \gamma'_\ep(A-),
\end{equation}
so $\gamma_\ep$ will be $C^1$ at $A$.

From property (c), \eqref{e5.104} and \eqref{e5.101}, we see that
\[
D_5 f(A+\ep) = \bar f^\ep_2(A+\ep) \int
_{A}^{A+\ep} \frac{-2\ep}{\rho
^2} \frac{1}{y^2(1-y)^2 (y-A)}
\frac{1}{\bar f^\ep_2(y)} \,dy,
\]
and this determines $D_5$ (and therefore $h_5$).

From (a) and \eqref{e5.103}, we find that
\[
D_4 = \int_B^{A+\ep}
h_5(y) \,dy,
\]
so that
%
%
\begin{equation}
\label{e5.107} \gamma_\ep(x) = \int_B^x
h_5(y) \,dy \qquad\mbox{for } A+\ep< x < B.
\end{equation}
From (b), \eqref{e5.107} and \eqref{e5.100}, we see that
\[
\int_B^{A+\ep} h_5(y) \,dy = \int
_A^{A+\ep} h_4(y) \,dy +
D_3,
\]
and this determines $D_3$.

Finally, from (d), \eqref{e5.98} and \eqref{e5.100}, we see that
\[
\frac{1}{\lambda} \ln(1-A) + D_1 = D_3,
\]
and this now determines $D_1$.

With the choice of constants $D_1,\dots,D_5$ above, 
we have determined a function $\gamma_\ep\dvtx[0,B] \to\IR$ which is
$C^1$ on $[0,B]$ 
and $C^2$ on $[0,A]$ and $[A, B]$.

We now turn to the study of $E_x(\tau^\ep)$. For $x \in[A,B]$, the
behavior of $p^\ep_t$ while $p^\ep_t \in[A,A+\ep[$ is somewhat
unusual, because of the square-root in the diffusion coefficient in the
s.d.e. \eqref{e5.76}. We are going to check below that in this
interval, $p^\ep_t - A$ is comparable to the time-change (under a
well-behaved time change) of a BESQ-process \cite{RY}, Chapter XI, so
it behaves essentially like the Cox--Ingersoll--Ross process \cite{LL},
Theorem~6.2.3, Proposition~6.2.4. In particular, $p^\ep_t \geq A$
since $x \geq A$, so $(p^\ep_t)$ never goes strictly below $A$ (though
it may hit $A$ and $A$ is instantaneously reflecting), and $A+\ep$ is
hit in finite time because $p^\ep_t$ is either recurrent or transient,
depending on the values of $\lambda$ and $\rho$. On the other hand,
$h^{(\ep)}(p^\ep_t) = 1$ while $p^\ep_t \in[A+\ep,B]$, so \eqref
{e5.76} simply describes there a diffusion with positive, bounded and
Lipschitz continuous drift and diffusion coefficients.

Regarding the behavior of $p^\ep_t$ while $p^\ep_t \in[A,A+\ep[$, set
$\tilde p^\ep_t = p^\ep_t - A$. By \eqref{e5.76},
%
%
\begin{equation}
\label{e5.48sept} d \tilde p^\ep_t = \lambda \bigl(1-A -
\tilde p^\ep_t \bigr) \,dt + \sqrt{\tilde
p^\ep_t} \sigma_t \,d\bar W^\ep_t,
\end{equation}
where
\[
\sigma_t = \rho\ep^{-1/2} \bigl(\tilde p^\ep_t
+A \bigr) \bigl(1-A - \tilde p^\ep_t \bigr).
\]
In particular, there are two positive and finite constants $c_\ep$ and
$C_\ep$ such that $c_\ep\leq\sigma_t \leq C_\ep$ as long as $p^\ep_t
\in[A,A+\ep[$. Define a martingale
%
%
\begin{equation}
\label{e5.49sept} M_t = \int_0^t
\sigma_s \,d\bar W^\ep_s,
\end{equation}
so that
\[
\langle M\rangle_t = \int_0^t
\sigma_s^2 \,ds.
\]
Define the increasing process $(\rho_t)$ so that $\langle M\rangle
_{\rho
_t} = t$, and notice that
\[
c_\ep^2 t \leq\langle M\rangle_t \leq
C_\ep^2t\quad \mbox{and}\quad \frac
{t}{C_\ep^2} t \leq
\rho_t \leq\frac{t}{c_\ep^2} t.
\]
By \eqref{e5.48sept} and \eqref{e5.49sept},
\[
\tilde p^\ep_t = \tilde p^\ep_0 +
\int_0^t \lambda \bigl(1-A - \tilde
p^\ep_s \bigr) \,ds + \int_0^t
\sqrt{\tilde p^\ep_s} \,dM_s.
\]
Using the time-change formulas for deterministic and stochastic
integrals (see Problem 4.5 and Proposition~4.8 in \cite{KS}, Chapter~3, we
see that
\[
\tilde p^\ep_t = \tilde p^\ep_0 +
\int_0^{\langle M\rangle_t} \lambda \bigl(1-A - \tilde
p^\ep_{\rho_s} \bigr) \frac{1}{\sigma^2_{\rho_s}} \,ds + \int
_0^{\langle M\rangle_t} \sqrt{q_s} \,d\bar
W^\ep_s.
\]
Setting $q_s = \tilde p^\ep_{\rho_s}$, so that $\tilde p^\ep_t =
q_{\langle M\rangle_t}$, we find that
\[
q_{\langle M\rangle_t} = \int_0^{\langle M\rangle_t} \lambda(1-A -
q_s) \frac{1}{\sigma^2_{\rho_s}} \,ds + \int_0^{\langle M\rangle_t}
\sqrt{q_s} \,d\bar W^\ep_s,
\]
and, setting $t=\rho_u$,
\[
q_u = \int_0^u \lambda(1-A -
q_s) \frac{1}{\sigma^2_{\rho_s}} \,ds + \int_0^u
\sqrt{q_s} \,d\bar W^\ep_s.
\]
The drift of $(q_u)$ is $\lambda(1-A - q_u)/\sigma_{\rho_u}^2 \geq
\lambda(1-A-\ep)/C_\ep^2$, so by the comparison theorem for s.d.e.'s
\cite{KS}, Chapter~5, Proposition~2.18, $q_u$ is greater than the
BESQ-process with drift $\lambda(1-A-\ep)/C_\ep^2$, hence $q_u \geq0$
a.s., or, equivalently, $p^\ep_t \geq A$ a.s.

We now apply It\^o's formula to $\gamma_\ep(p^\ep_{t\wedge\tau^\ep
})$, since $\gamma_\ep$ is $C^2$ on $[A,B]$,
\begin{eqnarray*}
\gamma_\ep \bigl(p^\ep_{t\wedge\tau^\ep} \bigr) &=&
\gamma_\ep \bigl(p^\ep_0 \bigr) + \int
_0^{t\wedge\tau^\ep} \gamma'_\ep
\bigl(p^\ep_s \bigr) \,dp^\ep_s +
\frac{1}{2}\int_0^{t\wedge\tau^\ep}
\gamma''_\ep \bigl(p^\ep_s
\bigr) \,d \bigl\langle p^\ep \bigr\rangle_s
\\
&=& \gamma_\ep \bigl(p^\ep_0 \bigr) + \int
_0^{t\wedge\tau^\ep} \gamma'_\ep
\bigl(p^\ep_s \bigr) \rho p^\ep_s
\bigl(1-p^\ep_s \bigr) \sqrt{h^{(\ep)}
\bigl(p^\ep_s \bigr)} \,d \bar W^\ep_s
\\
&&{} + \int_0^{t\wedge\tau^\ep} L\gamma_\ep
\bigl(p^\ep_s \bigr) \,ds.
\end{eqnarray*}
According to \eqref{e5.96}, $L\gamma_\ep(p^\ep_s) = -1$ for $s <
\tau
^\ep$, so, taking expectations, we find that
\[
E_x \bigl(\gamma_\ep \bigl(p^\ep_{t\wedge\tau^\ep}
\bigr) \bigr) = \gamma_\ep(x) - E_x \bigl(t\wedge
\tau^\ep \bigr),
\]
so
\[
E_x \bigl(t\wedge\tau^\ep \bigr) = - E_x
\bigl( \gamma_\ep \bigl(p^\ep_{t\wedge\tau^\ep} \bigr) \bigr)
+ \gamma_\ep(x).
\]

The right-hand side is bounded, so $\sup_{t \in\IR_+} E_x(t\wedge
\tau
^\ep) < +\infty$. By the monotone convergence theorem, $E_x(\tau^\ep
) <
+\infty$ as claimed [and in fact,\break  $E_x(\tau^\ep) = \gamma_\ep
(x)$], $x
\in[A,B]$.

For $x \in[0,A[$, we observe from \eqref{e5.76} that $p^\ep_t$ is
deterministic and increases at speed $\geq\lambda(1-A)$ until
reaching $A$. Thus $A$ is hit in less than some $\tau_0$ units of time,
and so
\[
E_x \bigl(\tau^\ep \bigr) \leq\tau_0 +
E_A \bigl(\tau^\ep \bigr) < +\infty.
\]
Finally, for $x \in[B,1]$, $\tau^\ep= 0$ $P_x$-a.s., so $E_x(\tau
^\ep
) = 0$. This proves Lemma~\ref{lem16}.
\end{pf}

%
\begin{lemma}\label{lem17} The function $g_\ep$ defined in
Proposition~\ref{lem15}
is the cost associated with the strategy $(h^\ep, \tau^\ep)$, that is,
for all $x \in[0,1]$, $g_\ep(x) = \tilde g_\ep(x) = E(C(h^\ep, \tau
^\ep
))$ [$\tilde g_\ep$ is defined in \eqref{e5.71}].
\end{lemma}

\begin{pf} According to Proposition~\ref{lem15}, $g_\ep$ satisfies the
assumptions of Lem\-ma~\ref{lem14}, and according to Lemma~\ref{lem16},
\eqref{e5.75} holds. Therefore, by Lem\-ma~\ref{lem14}, $g_\ep= \tilde
g_\ep$, and this proves Lemma~\ref{lem17}.
\end{pf}

\begin{pf*}{Proof of Theorem~\ref{thmvf1}} In case (a) of Theorem~\ref
{thmvf}, in view of the considerations that follow the proof of
Theorem~\ref{thmvf}, it remains only to prove that $\tilde g \leq g$.
By definition of $\tilde g$ and Lemma~\ref{lem17}, the inequality
$\tilde g \leq g_\ep$ holds. Since $g = \lim_{\ep\downarrow0} g_\ep$
by Proposition~\ref{lem15}, we conclude that $\tilde g \leq g$. This
completes the proof of Theorem~\ref{thmvf1} in case (a) of
Theorem~\ref{thmvf}.

The statement of Theorem~\ref{thmvf1} in case (b) of Theorem~\ref
{thmvf} follows from the fact that the optimal control system exhibited
in the proof of Theorem~\ref{thmvf} is (trivially) a strategy, which
then is necessarily optimal.
\end{pf*}

\section*{Acknowledgements} The authors thank two referees for a
very careful reading of the manuscript, which led to several important
clarifications.

%



\printaddresses

\begin{thebibliography}{24}

\bibitem{balmer1}
%
\begin{barticle}[mr]
\bauthor{\bsnm{Balmer},~\bfnm{D.~W.}\binits{D.~W.}}
(\byear{1975}).
\btitle{On a quickest detection problem with costly information}.
\bjournal{J. Appl. Probab.}
\bvolume{12}
\bpages{87--97}.
\bid{issn={0021-9002}, mr={0431582}}
\end{barticle}
%
\bptok{imsref}%
\endbibitem

\bibitem{balmer2}
%
\begin{barticle}[mr]
\bauthor{\bsnm{Balmer},~\bfnm{D.~W.}\binits{D.~W.}}
(\byear{1976}).
\btitle{On a quickest detection problem with variable monitoring}.
\bjournal{J. Appl. Probab.}
\bvolume{13}
\bpages{760--767}.
\bid{issn={0021-9002}, mr={0431583}}
\end{barticle}
%
\bptok{imsref}%
\endbibitem

\bibitem{BV1}
%
\begin{barticle}[mr]
\bauthor{\bsnm{Banerjee},~\bfnm{Taposh}\binits{T.}} \AND
\bauthor{\bsnm{Veeravalli},~\bfnm{Venugopal~V.}\binits{V.~V.}}
(\byear{2012}).
\btitle{Data-efficient quickest change detection with on-off
observation control}.
\bjournal{Sequential Anal.}
\bvolume{31}
\bpages{40--77}.
\bid{doi={10.1080/07474946.2012.651981}, issn={0747-4946}, mr={2899604}}
\end{barticle}
%
\bptok{imsref}%
\endbibitem

\bibitem{BV2}
%
\begin{barticle}[mr]
\bauthor{\bsnm{Banerjee},~\bfnm{Taposh}\binits{T.}} \AND
\bauthor{\bsnm{Veeravalli},~\bfnm{Venugopal~V.}\binits{V.~V.}}
(\byear{2013}).
\btitle{Data-efficient quickest change detection in minimax settings.}
\bjournal{IEEE Trans. Inform. Theory}
\bvolume{59}
\bpages{6917--931}.
\bid{mr={3106874}}
\end{barticle}
%
\bptok{imsref}%
\endbibitem

\bibitem{bather}
%
\begin{barticle}[auto]
\bauthor{\bsnm{Bather},~\bfnm{J.~A.}\binits{J.~A.}}
(\byear{1973}).
\btitle{An optimal stopping problem with costly information}.
\bjournal{Bull. Inst. Internat. Statist.}
\bvolume{45}
\bpages{9--24}.
\bid{mr={0402919}}
\end{barticle}
%
\bptok{imsref}%
\endbibitem

\bibitem{BK}
%
\begin{bmisc}[auto:STB|2014/02/12|14:17:21]
\bauthor{\bsnm{Bayraktar},~\bfnm{E.}\binits{E.}} \AND
\bauthor{\bsnm{Kravitz},~\bfnm{R.}\binits{R.}}
(\byear{2012}).
\bhowpublished{Quickest detection with discretely controlled
observations. Preprint.
Available at \arxivurl{arXiv:1212.4717v2}.}
\end{bmisc}
%
\bptok{imsref}%
\endbibitem

\bibitem{CD}
%
\begin{bbook}[mr]
\bauthor{\bsnm{Cairoli},~\bfnm{R.}\binits{R.}} \AND
\bauthor{\bsnm{Dalang},~\bfnm{Robert~C.}\binits{R.~C.}}
(\byear{1996}).
\btitle{Sequential Stochastic Optimization}.
\bpublisher{Wiley},
\blocation{New York}.
\bid{doi={10.1002/9781118164396}, mr={1369770}}
\end{bbook}
%
\bptok{imsref}%
\endbibitem

\bibitem{Chi}
%
\begin{barticle}[mr]
\bauthor{\bsnm{Chitashvili},~\bfnm{R.}\binits{R.}}
(\byear{1997}).
\btitle{On the nonexistence of a strong solution in the boundary
problem for a sticky {B}rownian motion}.
\bjournal{Proc. A. Razmadze Math. Inst.}
\bvolume{115}
\bpages{17--31}.
\bid{issn={1512-0007}, mr={1639096}}
\end{barticle}
%
\bptok{imsref}%
\endbibitem

\bibitem{dayanik}
%
\begin{barticle}[mr]
\bauthor{\bsnm{Dayanik},~\bfnm{Savas}\binits{S.}}
(\byear{2010}).
\btitle{Wiener disorder problem with observations at fixed discrete
time epochs}.
\bjournal{Math. Oper. Res.}
\bvolume{35}
\bpages{756--785}.
\bid{doi={10.1287/moor.1100.0471}, issn={0364-765X}, mr={2777513}}
\end{barticle}
%
\bptok{imsref}%
\endbibitem

\bibitem{durrett}
%
\begin{bbook}[mr]
\bauthor{\bsnm{Durrett},~\bfnm{Richard}\binits{R.}}
(\byear{1996}).
\btitle{Probability: Theory and Examples},
\bedition{2nd} ed.
\bpublisher{Duxbury Press},
\blocation{Belmont, CA}.
\bid{mr={1609153}}
\end{bbook}
%
\bptok{imsref}%
\endbibitem

\bibitem{elkaroui}
%
\begin{bincollection}[mr]
\bauthor{\bsnm{El Karoui},~\bfnm{N.}\binits{N.}}
(\byear{1981}).
\btitle{Les aspects probabilistes du contr\^ole stochastique}.
In \bbooktitle{Ninth {S}aint {F}lour {P}robability {S}ummer
{S}chool---1979 ({S}aint {F}lour, 1979)}.
\bseries{Lecture Notes in Math.}
\bvolume{876}
\bpages{73--238}.
\bpublisher{Springer},
\blocation{Berlin}.
\bid{mr={0637471}}
\end{bincollection}
%
\bptok{imsref}%
\endbibitem

\bibitem{GS}
%
\begin{bbook}[mr]
\bauthor{\bsnm{G{\={\i}}hman},~\bfnm{{\u{I}}.~{\=I}.}\binits{\u
{I}.~{\=I}.}} \AND
\bauthor{\bsnm{Skorohod},~\bfnm{A.~V.}\binits{A.~V.}}
(\byear{1972}).
\btitle{Stochastic Differential Equations}.
\bpublisher{Springer},
\blocation{Berlin}.
\bid{mr={0346904}}
\end{bbook}
%
\bptok{imsref}%
\endbibitem

\bibitem{IW}
%
\begin{bbook}[mr]
\bauthor{\bsnm{Ikeda},~\bfnm{Nobuyuki}\binits{N.}} \AND
\bauthor{\bsnm{Watanabe},~\bfnm{Shinzo}\binits{S.}}
(\byear{1981}).
\btitle{Stochastic Differential Equations and Diffusion Processes}.
\bseries{North-Holland Mathematical Library}
\bvolume{24}.
\bpublisher{North-Holland},
\blocation{Amsterdam}.
\bid{mr={0637061}}
\end{bbook}
%
\bptok{imsref}%
\endbibitem

\bibitem{KarShir}
%
\begin{barticle}[mr]
\bauthor{\bsnm{Karatzas},~\bfnm{Ioannis}\binits{I.}},
\bauthor{\bsnm{Shiryaev},~\bfnm{Albert~N.}\binits{A.~N.}} \AND
\bauthor{\bsnm{Shkolnikov},~\bfnm{Mykhaylo}\binits{M.}}
(\byear{2011}).
\btitle{On the one-sided {T}anaka equation with drift}.
\bjournal{Electron. Commun. Probab.}
\bvolume{16}
\bpages{664--677}.
\bid{doi={10.1214/ECP.v16-1665}, issn={1083-589X}, mr={2853104}}
\end{barticle}
%
\bptok{imsref}%
\endbibitem

\bibitem{KS}
%
\begin{bbook}[mr]
\bauthor{\bsnm{Karatzas},~\bfnm{Ioannis}\binits{I.}} \AND
\bauthor{\bsnm{Shreve},~\bfnm{Steven~E.}\binits{S.~E.}}
(\byear{1988}).
\btitle{Brownian Motion and Stochastic Calculus}.
\bseries{Graduate Texts in Mathematics}
\bvolume{113}.
\bpublisher{Springer},
\blocation{New York}.
\bid{doi={10.1007/978-1-4684-0302-2}, mr={0917065}}
\end{bbook}
%
\bptok{imsref}%
\endbibitem

\bibitem{LL}
%
\begin{bbook}[mr]
\bauthor{\bsnm{Lamberton},~\bfnm{Damien}\binits{D.}} \AND
\bauthor{\bsnm{Lapeyre},~\bfnm{Bernard}\binits{B.}}
(\byear{1996}).
\btitle{Introduction to Stochastic Calculus Applied to Finance}.
\bpublisher{Chapman \& Hall},
\blocation{London}.
\bid{mr={1422250}}
\end{bbook}
%
\bptok{imsref}%
\endbibitem

\bibitem{oksendal}
%
\begin{bbook}[mr]
\bauthor{\bsnm{{\O}ksendal},~\bfnm{Bernt}\binits{B.}}
(\byear{2003}).
\btitle{Stochastic Differential Equations: An Introduction with Applications},
\bedition{6th} ed.
\bpublisher{Springer},
\blocation{Berlin}.
\bid{doi={10.1007/978-3-642-14394-6}, mr={2001996}}
\end{bbook}
%
\bptok{imsref}%
\endbibitem

\bibitem{OS}
%
\begin{bbook}[mr]
\bauthor{\bsnm{{\O}ksendal},~\bfnm{Bernt}\binits{B.}} \AND
\bauthor{\bsnm{Sulem},~\bfnm{Agn{\`e}s}\binits{A.}}
(\byear{2007}).
\btitle{Applied Stochastic Control of Jump Diffusions},
\bedition{2nd} ed.
\bpublisher{Springer},
\blocation{Berlin}.
\bid{doi={10.1007/978-3-540-69826-5}, mr={2322248}}
\end{bbook}
%
\bptok{imsref}%
\endbibitem

\bibitem{PS}
%
\begin{bbook}[mr]
\bauthor{\bsnm{Peskir},~\bfnm{Goran}\binits{G.}} \AND
\bauthor{\bsnm{Shiryaev},~\bfnm{Albert}\binits{A.}}
(\byear{2006}).
\btitle{Optimal Stopping and Free-Boundary Problems}.
\bpublisher{Birkh\"auser},
\blocation{Basel}.
\bid{mr={2256030}}
\end{bbook}
%
\bptok{imsref}%
\endbibitem

\bibitem{PS2}
%
\begin{barticle}[mr]
\bauthor{\bsnm{Peskir},~\bfnm{G.}\binits{G.}} \AND
\bauthor{\bsnm{Shiryaev},~\bfnm{A.~N.}\binits{A.~N.}}
(\byear{2000}).
\btitle{Sequential testing problems for {P}oisson processes}.
\bjournal{Ann. Statist.}
\bvolume{28}
\bpages{837--859}.
\bid{doi={10.1214/aos/1015952000}, issn={0090-5364}, mr={1792789}}
\end{barticle}
%
\bptok{imsref}%
\endbibitem

\bibitem{protter}
%
\begin{bbook}[mr]
\bauthor{\bsnm{Protter},~\bfnm{Philip~E.}\binits{P.~E.}}
(\byear{2004}).
\btitle{Stochastic Integration and Differential Equations},
\bedition{2nd} ed.
\bseries{Applications of Mathematics (New York)}
\bvolume{21}.
\bpublisher{Springer},
\blocation{Berlin}.
\bid{mr={2020294}}
\end{bbook}
%
\bptok{imsref}%
\endbibitem

\bibitem{RY}
%
\begin{bbook}[mr]
\bauthor{\bsnm{Revuz},~\bfnm{Daniel}\binits{D.}} \AND
\bauthor{\bsnm{Yor},~\bfnm{Marc}\binits{M.}}
(\byear{1991}).
\btitle{Continuous Martingales and {B}rownian Motion}.
\bseries{Grundlehren der Mathematischen Wissenschaften}
\bvolume{293}.
\bpublisher{Springer},
\blocation{Berlin}.
\bid{doi={10.1007/978-3-662-21726-9}, mr={1083357}}
\end{bbook}
%
\bptok{imsref}%
\endbibitem

\bibitem{shiryaev}
%
\begin{bbook}[mr]
\bauthor{\bsnm{Shiryayev},~\bfnm{A.~N.}\binits{A.~N.}}
(\byear{1978}).
\btitle{Optimal Stopping Rules}.
\bpublisher{Springer},
\blocation{Berlin}.
\bid{mr={0468067}}
\end{bbook}
%
\bptok{imsref}%
\endbibitem

\bibitem{warren}
%
\begin{bincollection}[mr]
\bauthor{\bsnm{Warren},~\bfnm{J.}\binits{J.}}
(\byear{1999}).
\btitle{On the joining of sticky {B}rownian motion}.
In \bbooktitle{S\'eminaire de {P}robabilit\'es, {XXXIII}}.
\bseries{Lecture Notes in Math.}
\bvolume{1709}
\bpages{257--266}.
\bpublisher{Springer},
\blocation{Berlin}.
\bid{doi={10.1007/BFb0096515}, mr={1767999}}
\end{bincollection}
%
\bptok{imsref}%
\endbibitem
\end{thebibliography}
\end{document}